\documentclass[12pt]{amsart}
\usepackage{amsmath,amssymb}

\newtheorem{thm}{Theorem}[section]
\newtheorem{cor}[thm]{Corollary}
\newtheorem{lemma}[thm]{Lemma}
\newtheorem{prop}[thm]{Proposition}

\theoremstyle{definition}
\newtheorem{definition}[thm]{Definition}

\newtheorem{case}{Case}

\theoremstyle{remark}

\newtheorem{ques}[thm]{Question}

\numberwithin{equation}{section}

\newcommand{\thmref}[1]{Theorem~\ref{#1}}
\newcommand{\secref}[1]{\S\ref{#1}}
\newcommand{\lemref}[1]{Lemma~\ref{#1}}
\newcommand{\propref}[1]{Prop\-o\-si\-tion~\ref{#1}}
\newcommand{\corref}[1]{Cor\-ol\-lary~\ref{#1}}

\newcommand{\quesref}[1]{Question~\ref{#1}}
\newcommand{\defref}[1]{Definition~\ref{#1}}

\makeatletter
\newif\ifdraft
\draftfalse
\newcount\hours
\newcount\minutes       %   For computing the time of day on
\def\timeofday{%    Must be computed when called if preloaded
\hours=\time
\minutes=\hours
\divide\hours by60
\multiply\hours by60
\advance\minutes by-\hours
\divide\hours by60
\ifnum\hours>9\else0\fi\the\hours:\ifnum\minutes>9\else
0\fi\the\minutes}

\def\calendar{\ifcase\month\or January\or February\or March\or April\or
May\or June\or July\or August\or September\or October\or November\or
December\fi\space\number\day, \number\year}
\def\today{\calendar\quad \timeofday}

\def\T{\le_T}
\def\m{\le_m}
\def\eT{\equiv_T}

\newcommand{\jump}[2]{{#1}^{#2}}
\newcommand{\td}[1]{\mathbf{#1}}
\newcommand{\red}[3]{\{#1\}^{#2}_{#3}}
\def\zerop{0^\prime}
\def\zeropp{0^{\prime \prime}}
\newcommand{\homo}[1]{H({#1})}
\newcommand{\ktup}[2]{[{#1}]^{#2}}

\newcommand{\fin}[1]{\mathrm{Fin}^{#1}}
\newcommand{\low}[1]{\mathrm{low}_{#1}}
\def\join{\oplus}
\newcommand{\color}[1]{\mathcal{#1}}

\begin{document}

%\topmatter
\title[Generalized cohesiveness]{Generalized cohesiveness}

%author one info
\author{Tamara Hummel}
\address{Department of Mathematics\\
         Allegheny College\\
         520 N.~Main St.\\
         Meadville, PA  16335}
%For AMS:  where research conducted.  Divide into lines using \\
%\curaddress{}
\email{thummel@alleg.edu}
%\thanks{}

%author two info
\author{Carl G. Jockusch, Jr.}
\address{Department of Mathematics\\
         University of Illinois\\
         1409 W.~Green St.\\
         Urbana, IL  61801}
\email{jockusch@math.uiuc.edu}
\thanks{Some of the material in this paper originally appeared as part of
the first author's Ph.D.~thesis.  Jockusch was partially supported by NSF
Grant DMS 95-03398}


\date{\today}

\dedicatory{}
%\endtopmatter

\begin{abstract}
We study some generalized notions of cohesiveness which arise
naturally in connection with effective versions of Ramsey's Theorem.
An infinite set $A$ of natural numbers is $n$--cohesive (respectively,
$n$--r--cohesive) if $A$ is almost homogeneous for every computably
enumerable (respectively, computable) $2$--coloring of the
$n$--element sets of natural numbers.  (Thus the $1$--cohesive and
$1$--r--cohesive sets coincide with the cohesive and r--cohesive sets,
respectively.)  We consider the degrees of unsolvability and
arithmetical definability levels of $n$--cohesive and $n$--r--cohesive
sets.  For example, we show that for all $n \ge 2$, there exists a
$\Delta^0_{n+1}$ $n$--cohesive set.  We improve this result for $n =
2$ by showing that there is a $\Pi^0_2$ $2$--cohesive set.  We show
that the $n$--cohesive and $n$--r--cohesive degrees together form a
linear, non--collapsing hierarchy of degrees for $n
\geq 2$.  In addition, for $n \geq 2$ we characterize the jumps of
$n$--cohesive degrees as exactly the degrees ${\bf \geq
\jump{0}{(n+1)}}$ and show that each $n$--r--cohesive degree has jump
${\bf > \jump{0}{(n)}}$.
\end{abstract}

\maketitle

\section{Introduction}

We study a hierarchy of generalized notions of cohesiveness,
which arises naturally in connection with
effective versions of Ramsey's Theorem.  For a set $X \subseteq
\omega$, let $\ktup{X}{n}$ denote the class of all $n$--element
subsets of  $X$.  A {\em $k$--coloring\/} $\color{C}$ of
$\ktup{X}{n}$ is a function $\color{C} : \ktup{X}{n} \to \{0, 1,
\dots, k-1\}$; $n$ is called the {\em  exponent\/} of the
coloring.  A set $A \subseteq \omega$ is {\em  homogeneous\/} for
a coloring $\color{C}$ of exponent $n$ if $\color{C}$ is constant
on $\ktup{A}{n}$; i.e., all $n$--element subsets of $A$ are
assigned the same color by  $\color{C}$.  The infinite form of
Ramsey's Theorem \cite{Ra} (Theorem~A) states that for any
infinite set $X$  and any $k$--coloring $\color{C}$ of
$\ktup{X}{n}$, there is an infinite set  $A \subseteq X$ such
that $A$ is homogeneous for $\color{C}$.

Effective versions of Ramsey's Theorem were studied in detail by Jockusch
in \cite{Jo1}, which considered the degrees of
unsolvability and arithmetical  definability properties of
infinite homogeneous sets for computable (recursive)  colorings.  In this
paper, we study both computable and computably  enumerable (c.e.) $2$--colorings 
of $\ktup{\omega}{n}$.

\begin{definition}  A $2$--coloring $\color{C}$ of $\ktup{\omega}{n}$ is {\em  
computably enumerable (c.e.)\/} if either $\{D \in \ktup{\omega}{n} \mid
\color{C}(D) = 0\}$ or $\{D \in \ktup{\omega}{n} \mid \color{C}(D) = 1\}$ is  
computably enumerable, under a suitable encoding of $\ktup{\omega}{n}$.

We often identify a c.e.~set $\color{C} \subseteq
\ktup{\omega}{n}$ with  its characteristic function and so
consider it to be a c.e.~$2$--coloring of  $\ktup{\omega}{n}$.
\end{definition}

Recall that an infinite set $A \subseteq \omega$ is {\em cohesive\/} if for  
all c.e.~sets $W$, $A \subseteq^* W$ or $A \subseteq^* \overline{W}$, where  
$X \subseteq^* Y$ means there exists a finite set $F$ such that $X - F
\subseteq Y$.  Similarly, an infinite set $A$ is {\em r--cohesive\/} if for  
all computable sets $R$, $A \subseteq^* R$ or $A \subseteq^*
\overline{R}$.  Using the language of Ramsey's Theorem, an infinite
set $A$ is cohesive if and only if it is ``almost homogeneous'' for
every c.e.~$2$--coloring of $\ktup{\omega}{1}$.

\begin{definition}  A set $A$ is {\em almost homogeneous\/} for a coloring
$\color{C}$ if there is a finite set $F$ such that $A - F$ is homogeneous for  
$\color{C}$.
\end{definition}

By considering $2$--colorings of $\ktup{\omega}{n}$, we generalize the notion  
of cohesiveness to a natural hierarchy of strong forms of cohesiveness.

\begin{definition}  An infinite set $A$ is {\em $n$--cohesive\/}
(respectively, {\em $n$--r--cohesive\/}) if $A$ is almost homogeneous for every  
c.e.~(respectively, computable) $2$--coloring of $\ktup{\omega}{n}$.
\end{definition}

Clearly, the $1$--cohesive sets are precisely the cohesive sets,
and the $1$--r--cohesive sets are exactly the r--cohesive sets.
It is easily seen that every $(n+1)$--r--cohesive set is $n$--cohesive,
and it is obvious by definition that every $n$--cohesive set is
$n$--r--cohesive.  It will follow from  \corref{cor:stronghier} (and the existence
of a set which is r--cohesive but not cohesive) that the converse
of each of these implications fails for each $n \geq 1$.   Thus we are
studying a linearly ordered proper hierarchy of ever stronger versions
of cohesiveness.

In this paper, we consider the degrees of unsolvability and arithmetical
definability properties of $n$--cohesive and $n$--r--cohesive sets.  Our
notation is consistent with that of Soare \cite{So1}, except that we
modify traditional terminology as suggested by Soare \cite{So2}.  In
particular, sets and functions traditionally called ``recursive'' are
here called ``computable,'' and sets traditionally called ``recursively
enumerable'' (or ``r.e.'') are here called ``computably enumerable''
(or ``c.e.'').

\section{Effective versions of Ramsey's Theorem}

In this section, we consider some of Jockusch's results \cite{Jo1} concerning  
effective versions of Ramsey's Theorem, as well as some generalizations of
those results, which will be needed in the sequel.  The following result
precisely locates in the arithmetical hierarchy which infinite homogeneous
sets are guaranteed to exist for a computable $2$--coloring.

\begin{thm}[Jockusch]\label{thm:basicjockusch}
\begin{enumerate}
\item  (\cite{Jo1}, Theorems~4.2 and 5.5) For all $n$ and  $k$, and for any  
computable $k$--coloring of
$\ktup{\omega}{n}$, there exists an infinite $\Pi^0_n$ homogeneous set.
\item  (\cite{Jo1}, Theorems~3.1 and 5.1)  For all $n \ge 2$, there exists a  
computable $2$--coloring of
$\ktup{\omega}{n}$ which has no infinite $\Sigma^0_n$ homogeneous set.
\end{enumerate}
\end{thm}

The proof of the first part of \thmref{thm:basicjockusch} for $n=2$, which
is a finite injury priority argument with a $\zerop$ oracle, is easily
modified for c.e.~$2$--colorings of $\ktup{\omega}{2}$.

\begin{thm}\label{thm:pi2re}  Every c.e.~$2$--coloring of $\ktup{\omega}{2}$  
has an infinite $\Pi^0_2$ homogeneous set.
\end{thm}

\begin{proof}[Proof Sketch]  Let a c.e.~$2$--coloring $\color{C}$ of
$\ktup{\omega}{2}$ be given as a red--blue coloring, where the set of all red  
pairs of $\color{C}$ is c.e., and the set of all blue pairs of $\color{C}$ is  
co--c.e.

The proof is a minor adjustment of Jockusch's proof in \cite{Jo1}.  The idea  
is to define an increasing sequence of numbers $\{a_n\}_{n \in \omega}$ and  
a red--blue coloring of the $a_n$'s such that for all $i < j$, the color of  
$\{a_i,a_j\}$ in the given c.e.~$2$--coloring is the same as the color of
$a_i$.  To make the set $\Pi^0_2$, it is initially assumed that each $a_i$
can be colored red.    If this assumption is later found to be incorrect,
because, for example, only finitely many numbers make a red pair with $a_i$,  
then the color of $a_i$ is changed from red to blue, and the part of the
sequence
which is constructed based on the incorrect color of $a_i$ is destroyed.  The  
construction is a movable marker construction using a $\zerop$ oracle; let
$a_i^s$ denote the position of marker $\Lambda_i$ at the beginning of stage  
$s$.

Jockusch's original proof used the notion of a $k$--acceptable number; in his  
original proof, a number $c$ is $k$--acceptable at stage $s$ if for all $i <  
k$, $a_i^s$ is defined, $a_i^s < c$, and the color of $\{a_i^s, c\}$ (in the  
given {\em computable\/} $2$--coloring) is the same as the color of $a_i^s$.    
To ensure that the construction requires only a $\zerop$ oracle, the notion  
of $k$--acceptability is reworded in this proof to compensate for
c.e.~$2$--colorings; a number $c$ is {\em $k$--acceptable at stage $s$} if for  
all $i<k$, $a_i^s$ is defined, $a_i^s < c$, and $\{a_i^s,c\}$ is red if
$a_i^s$ is red.

In addition, say that a number $c$ is {\em free} at $s$ if it has never been  
the position of a marker prior to stage $s$, and $c \ge s$.  Note that when  
$s$ is fixed, to say ``$c$ is free and $k$--acceptable at $s$'' is a
$\Sigma^0_1$ predicate.

\noindent {\em Construction.}

\noindent {\em Stage $s \ge 0$.}  Inductively assume that there exists a
number $n(s)$ such that the markers currently having a position are exactly  
the $\Lambda_i$, for $i< n(s)$.

\begin{case}  There exists a number which is free and $n(s)$--acceptable at $s$.
\end{case}

 Attach the marker $\Lambda_{n(s)}$ to the least such number $c$ and color $c$ red.

\begin{case}  Otherwise.  (Correct a mistake.)
\end{case}

Let $j(s)$ be the largest number $j$ such that there exists a number which
is free and $j$--acceptable at $s$.  Note that such a number $j$ exists
because every number is $0$--acceptable at stage $s$, and that $j(s)<n(s)$.   
Change the color of $a_{j(s)}^s$ and detach all markers $\Lambda_i$ for $j(s)  
< i < n(s)$.

The construction requires only a $\zerop$ oracle, as the noncomputable
questions in it ask whether certain given $\Sigma^0_1$ sets are nonempty.     
The new notion of $k$--acceptability suffices to ensure that the construction  
succeeds, as the only way that the color of $a_i^s$ can be changed from red  
to blue at stage $s$ is if there exist numbers which are free and
$i$--acceptable at $s$, but none of these numbers makes a red pair with
$a_i^s$.  Hence, all of these numbers make a blue pair with $a_i^s$.  Since  
the color of $a_i^s$ can never change back to red (although the marker
$\Lambda_i$ may later be detached from $a_i^s$), the definition of ``$c$ is  
$k$--acceptable at $s$'' has the property that if $i < k$ and $a_i^s$ is blue,  
then $\{a_i^s,c\}$ is blue.

The proofs of the following lemmas go through as in \cite{Jo1}.

\begin{lemma}  For all $k$, $\lim_s a_k^s = a_k$ exists, and the color of
$a_k$ can change from red to blue only (hence the color of $a_k$ stabilizes).
\end{lemma}

\begin{lemma}  If $i<j$, then the pair $\{a_i,a_j\}$ has the same color as
the eventual color of $a_i$.
\end{lemma}

Then define $M = \{a_i \mid i \in \omega\}$, $R = \{a_i \in M \mid a_i
\mbox{ is eventually red}\}$, and $B = \{a_i \in M \mid a_i \mbox{ is
eventually blue}\}$.  The set $M$ is infinite since $\{a_i\}_{i \in \omega}$  
is an increasing sequence.  The sets $R$ and $B$ are each homogeneous for the  
$2$--coloring $\color{C}$.  As in \cite{Jo1}, the sets $M$ and $R$ are each  
$\Pi^0_2$, so that if $R$ is infinite, then $R$ is the desired infinite
$\Pi^0_2$ homogeneous set.  If $R$ is finite, then $B = M-R$ is $\Pi^0_2$ and  
is the desired infinite homogeneous set.
\end{proof}

This result naturally leads to the following question.

\begin{ques}\label{ques:biggie}  For which $n$ does every
c.e.~$2$--coloring of  $\ktup{\omega}{n}$ have an infinite
$\Pi^0_n$ homogeneous set? \end{ques}

\quesref{ques:biggie} remains open;
in fact, it is unknown whether every  c.e.~$2$--coloring of
$\ktup{\omega}{3}$ has an infinite $\Pi^0_3$ homogeneous  set.
However, we improve the least known arithmetical complexity of
infinite homogeneous sets for c.e.~$2$--colorings of
$\ktup{\omega}{n}$, $n  \ge 3$, in \secref{sec:general}.

Another result of Jockusch considers the degree of homogeneous sets.

\begin{thm}[Jockusch, \cite{Jo1} (Corollary~4.7)]\label{thm:jump}  Every
computable
$2$--coloring of $\ktup{\omega}{2}$ has an infinite homogeneous set $A$ such  
that $\jump{A}{\prime} \T \zeropp$.
\end{thm}

This theorem was extended by Hummel \cite{Hu1} (Corollary~4.16) to
c.e.~$2$--colorings of $\ktup{\omega}{2}$ using a non--uniform
argument quite different from the uniform proof of \thmref{thm:jump}
given in \cite{Jo1}.  This extension will appear in a later
publication.  Our next topic is the question of whether the
case $n = k = 2$ of part 1 of \thmref{thm:basicjockusch} and
\thmref{thm:jump} can be combined.

\begin{ques}\label{ques:combine}  Does every computable $2$--coloring of
$\ktup{\omega}{2}$ have an infinite $\Pi^0_2$ homogeneous set $A$ with
$\jump{A}{\prime} \T \zeropp$?
\end{ques}

\subsection{Effective $\Delta^0_1$--immunity}

We give a partial answer to \quesref{ques:combine} by proving a theorem that  
extends the second part of \thmref{thm:basicjockusch}, which states that
there exists a computable $2$--coloring of $\ktup{\omega}{2}$ which has no
infinite $\zerop$--computable homogeneous set.  Jockusch's proof of this
result uses a construction which  yields a computable coloring, all of whose  
infinite homogeneous sets have a special ``effective immunity''
property.  When this  is combined with Martin's theorem
\cite{Ma1} that effectively simple sets are complete, we get a
partial negative answer to \quesref{ques:combine}, as well as a
result which can be directly applied to $2$--cohesive sets.

First recall that an infinite set $A \subseteq \omega$ is {\em immune\/} if  
for each c.e.~set $W_e$, if $W_e \subseteq A$, then $W_e$ is finite.
``Effectivizing'' this definition, we say that an infinite set $A$ is {\em
effectively immune\/} if there exists a (total) computable function $f$ such  
that for all $e$, if $W_e \subseteq A$, then $|W_e| \le f(e)$.  We define a  
new notion of effective immunity where, instead of considering c.e.~subsets  
of $A$, we consider computable subsets of $A$ given by $\Delta^0_1$--indices.   
In what follows, the notation $\downarrow$ abbreviates the phrase ``is
defined,'' while $\uparrow$ abbreviates the phrase ``is undefined.''

\begin{definition}
An infinite set $A$ is {\em effectively $\Delta^0_1$--immune\/} if
there exists a computable partial function $\psi$ such that for
all $a$ and $b$, \[(W_a \subseteq A \land \overline{W_a} = W_b)
\implies \psi(\langle a, b \rangle)\downarrow \land \ |W_a| \le
\psi(\langle a, b \rangle ).\] \end{definition}

At first glance, this notion seems quite different from that of
effective immunity.  Intuitively, the property that $\psi$ above is
partial seems necessary, as there is no reason to expect such a
function to be defined on pairs $\langle a, b \rangle$ which are not
$\Delta^0_1$--indices of computable sets.  However, it turns out that
the notion of effective $\Delta^0_1$--immunity is precisely the same
as that of effective immunity.

\begin{lemma}\label{lem:effimmune}
The following are equivalent:
\begin{enumerate}
\item $A$ is effectively $\Delta^0_1$--immune via a computable
partial function.
\item $A$ is effectively $\Delta^0_1$--immune
via a computable total function.
\item $A$ is effectively immune.
\end{enumerate}
\end{lemma}

\begin{proof} ((1) $\implies$ (2))  \ Let $A$ be effectively
$\Delta^0_1$--immune via the computable partial function $\psi$.  We define a  
total computable function $f$ with the desired properties.  Given $\langle a,  
b \rangle$, to define $f(\langle a, b \rangle )$ we construct auxiliary
c.e.~sets $W_c$ and $W_d$.  By the double recursion theorem (see \cite{So1},  
Exercise~II.3.15), we can give ourselves their indices $c$ and $d$ in
advance, and use them effectively in the construction.  To compute $f(\langle  
a, b \rangle)$, enumerate the numbers $0, 1, 2, \dots$ into $W_d$, one
number per stage, until the first stage $s_0$ occurs such that
$\psi_{s_0}(\langle c, d  \rangle )\downarrow$.  Define
$f(\langle a, b \rangle )$ to be $\psi_{s_0}(\langle c, d
\rangle) + |W_{d, s_0}|$.  Continue the construction of $W_c$ and
$W_d$ by  then enumerating all elements of $W_a$ not already in
$W_{d, s_0}$ into $W_c$  and all elements of $W_b$ into $W_d$.
The partial function $f$ is computable since the double recursion
theorem holds effectively; i.e., $c$ and $d$ may be
effectively computed from $a$ and $b$. (See the double recursion
theorem with parameters in \cite{So1},  Exercise~II.3.15).)

Note that $f$ is total, since if $f(\langle a, b
\rangle)\uparrow$, then  $W_c = \emptyset$, and $W_d =
\omega$.   Hence $W_c = \overline{W_d} \subseteq A$, which
implies that $\psi(\langle  c, d \rangle)\downarrow$, and so
$f(\langle a, b \rangle)\downarrow$, a  contradiction.

Next, suppose that $W_a = \overline{W_b} \subseteq A$.  We show
$|W_a| \le  f(\langle a, b \rangle)$.  First note that
$\psi(\langle a, b  \rangle)\downarrow$ and $|W_a| \le
\psi(\langle a, b \rangle)$.  Also, $f(\langle a, b \rangle)
= \psi(\langle c,  d \rangle) + |W_{d, s_0}|$, where $s_0$ is the
stage at which $f(\langle a, b  \rangle)$ becomes defined.  Then
$W_c \subseteq A$, since $W_c \subseteq  W_a$.  Also, $W_c \cap
W_d = \emptyset$, since $W_a \cap W_b = \emptyset$ and  only
elements in $W_a$ not already in $W_{d,s_0}$ are enumerated into
$W_c$.
 To see that $W_c \cup W_d = \omega$, note that $W_a \cup W_b = \omega$,
every element $x \in W_b$ is enumerated into  $W_d$, and every element of
$W_a$ is enumerated into either $W_{d,s_0}$ or $W_c$.  Thus we have $W_c =
\overline{W_d} \subseteq A$, and so $|W_c| \le \psi(\langle c, d \rangle)$.   
Finally, $W_a \subseteq W_c \cup W_{d,s_0}$, so
\begin{align*}
|W_a| &\le |W_c| + |W_{d,s_0}|\\
&\le \psi(\langle c, d \rangle) + |W_{d,s_0}|\\
&= f(\langle a, b \rangle)
\end{align*}
as desired.

((2) $\implies$ (3))  \ Let $A$ be effectively $\Delta^0_1$--immune via the  
total function $f$.  We define a computable function $g$
which witnesses that $A$ is effectively immune.  Given $a$, we
define $g(a) = f(c,d)$, where $W_c$ and $W_d$ are auxiliary c.e.~sets
to be defined below.  As in the proof of ((1) $\implies$ (2))
above, by the double recursion theorem we may use $c$ and $d$
effectively in the definition of $W_c$ and $W_d$.   To
enumerate $W_c$ and $W_d$, compute $f(c,d)$  and search for a
stage $s_0$ such that $|W_{a,s_0}| > f(c,d)$.   When such a stage
$s_0$ is found (if ever), put all elements of $W_{a,s_0}$ into
$W_c$ and all other elements of $\omega$ into $W_d$.   The
function $g$ is total since $f$ is, and $g$ is computable by the
uniformity of the double recursion theorem.  To show that $g$
witnesses the effective immunity of $A$, assume for a
contradiction that $W_a \subseteq A$ and $|W_a| > g(a) = f(c,d)$.
Then $s_0$ exists, and $W_c \subseteq W_a \subseteq A$.  Further,
it is clear that $W_d = \overline{W_c}$.  Hence, since $f$
witnesses that $A$ is effectively $\Delta^0_1$--immune, $|W_c| <
f(c,d)$.  But $|W_c| = |W_{a, s_0}| > f(c,d)$ by choice of $s_0$,
so we have a contradiction.

((3) $\implies$ (1)) is clear.

\end{proof}

It was shown by Jockusch in Theorem~3.1 of \cite{Jo1} that there
is a computable $2$--coloring $\color{C}$ of  $\ktup{\omega}{2}$
with no infinite homogeneous set $A \T \zerop$.  By analyzing
that construction and applying the theorem just proved, we now
obtain a stronger result.

\begin{thm}\label{thm:alleffimmunehomos}  There exists a computable
$2$--coloring $\color{C}$ of $\ktup{\omega}{2}$ such that every infinite
homogeneous set $A$ is effectively immune relative to $\zerop$.
\end{thm}

\begin{proof}  By the proof of the Shoenfield Limit Lemma (see
\cite{So1}, Lemma~III.3.3),  there exists a uniformly computable
sequence of $\{0, 1\}$--valued computable  functions $\{f_e\}_{e
\in \omega}$ such that every set $A \T \zerop$ is $A =  \lim f_e$
for some $e$; that is, for all $x$, $A(x) = \lim_s f_e(x,s)$.  As
in \cite{Jo1}, we write $A_e = \lim f_e$, where $A_e$ is
undefined if for  some $x$, $\lim_s f_e(x,s)\uparrow$.  We let
$\homo{\color{C}}$ denote the  class of all infinite homogeneous
sets for a coloring $\color{C}$.  The construction used to prove
Theorem~3.1 of \cite{Jo1} shows that there is a computable
$2$--coloring $\color{C}$ of  $\ktup{\omega}{2}$ such that, for all $A$,
\[(A_e\downarrow \land \ A_e \subseteq A \in \homo{\color{C}})
\implies  |A_e| \le 2e+1.\]

To see that this suffices to prove the theorem, note
that there exists a total computable function $f$ such that if $\langle a, b  
\rangle$ is a $\Delta^0_1$--index, relative to $\zerop$, of a
$\zerop$--computable set $B$, then $B = A_{f(\langle a, b \rangle)}$.  The
construction shows that if $B \subseteq A \in \homo{\color{C}}$, then $|B|
\le 2 f(\langle a, b \rangle) + 1$; i.e., all $A \in \homo{\color{C}}$ are
effectively $\Delta^0_1$--immune, relative to $\zerop$, via the function
$\psi(\langle a, b \rangle) = 2 f(\langle a, b \rangle) + 1$.
Hence all such $A$ are effectively immune relative to $\zerop$ by
\lemref{lem:effimmune}, relativized to $\zerop$.
\end{proof}

We thus have our partial negative answer to \quesref{ques:combine}.

\begin{cor}\label{cor:nopi2andjump}  There exists a computable $2$--coloring  
of $\ktup{\omega}{2}$ such that for all infinite homogeneous sets $A$, if $A$  
is $\Pi^0_2$, then $\zeropp \T A \join \zerop$.
\end{cor}

\begin{proof}  Take the $2$--coloring constructed in the theorem.  A theorem  
of Martin \cite{Ma1} states that if a set $A$ is $\Pi^0_1$ and effectively
immune, then $\zerop \T A$.  The corollary then follows immediately from
\thmref{thm:alleffimmunehomos} and Martin's theorem, both
relativized to $\zerop$. \end{proof}

\subsection{Stable colorings}

In what follows, we will require the notion of a stable coloring.

\begin{definition}\label{def:stable}  Let $A$ be an infinite subset of
$\omega$.  A
$k$--coloring $\color{C}$ of $\ktup{A}{n+1}$ is {\em stable\/} if
for all $D \in \ktup{A}{n}$ there exists $a_0
\in \omega$ such that for all $a \ge a_0$ with $a \in A$, $\color{C}(D
\cup \{a\}) = \color{C}(D \cup \{a_0\})$.
\end{definition}

Intuitively, a $k$--coloring $\color{C}$ of $\ktup{\omega}{n+1}$ is stable if  
the color of an unordered $(n+1)$--tuple depends ultimately only on the least  
$n$ elements of the tuple.  Stable partitions play a
crucial role in the proof of Cholak, Jockusch, and Slaman
\cite{ChJoSl} that every computable $2$--coloring of pairs has
an infinite $\low{2}$ homogeneous set.  We give here some
easy results about stable colorings.  Note that a stable
coloring of $\ktup{\omega}{2}$ tends to be  better behaved than
an arbitrary coloring; it is unknown whether a stable coloring
of $\ktup{\omega}{n}$ is better behaved than an arbitrary
coloring  when $n \ge 3$.  As we have already remarked, every
c.e.~$2$--coloring of  $\ktup{\omega}{2}$ has an infinite
$\Pi^0_2$ homogeneous set, and this result  is best possible by
\thmref{thm:basicjockusch}.  However, by imposing  stability on
$2$--colorings, we can improve the arithmetical complexity of the
homogeneous set.

\begin{prop}\label{prop:stabled2}  Every c.e.~stable $2$--coloring of
$\ktup{\omega}{2}$ has an infinite $\Delta^0_2$ homogeneous set.
\end{prop}

\begin{proof}  Let a c.e.~stable $2$--coloring $\color{C}$ of
$\ktup{\omega}{2}$ be given as a red--blue coloring of pairs, where the set of  
red pairs of $\color{C}$ is c.e.~and the set of blue pairs of $\color{C}$ is  
co--c.e.  Consider the $\Sigma^0_2$ set
\[A = \{a \mid (\exists b_0)(\forall b \ge b_0)[\{a,b\} \mbox{ is blue}]\}.\] 

\setcounter{case}{0}
\begin{case}  $A$ is finite.
\end{case}
We show in this case that there exists an infinite computable red homogeneous  
set.  Let $a > \max(A)$.  By definition of $A$, there exist infinitely many  
numbers $b > a$ such that $\{a,b\}$ is red, and hence by stability of
$\color{C}$, we have that for all $a > \max(A)$, $\{a,b\}$ is red for
sufficiently large $b$.  It is now easy to construct an infinite computable  
red homogeneous set.

\begin{case}  $A$ is infinite.
\end{case}
We first show that if $C \subseteq A$ is infinite, then there exists an
infinite blue homogeneous set $D \T C \join \zerop$.  Let $d_0 = (\mu x)[x
\in C]$.  Since $d_0 \in A$, $\{d_0, b\}$ is blue for sufficiently large
numbers $b$.  Search $\zerop$--effectively for a number $b_0$ such that for  
all $b \ge b_0$, $\{d_0,b\}$ is blue.  We let $d_1 = (\mu x)[x \in C \land x  
> \max\{b_0,d_0\}]$.  Next, we search $\zerop$--effectively for a number $b_1$  
such that for all $b \ge b_1$, $\{d_0,b\}$ and $\{d_1,b\}$ are blue.  We let  
$d_2 = (\mu x)[x \in C \land x > \max\{b_1,d_1\}]$.  We can continue in this  
fashion to enumerate an infinite set $D = \{d_0, d_1, \dots\}$ in increasing  
order, using oracles for $C$ and $\zerop$.

Next, since $A$ is an infinite $\Sigma^0_2$ set, there exists an infinite
$\Delta^0_2$ subset $B \subseteq A$.  By the above, it follows that there
exists an infinite blue homogeneous set $D \T B \join \zerop \T \zerop$.
\end{proof}

The next result shows that \propref{prop:stabled2} is best possible for
c.e.~stable $2$--colorings of $\ktup{\omega}{2}$.

\begin{prop}  There exists a computable stable $2$--coloring of
$\ktup{\omega}{2}$ with no infinite $\Pi^0_1$ homogeneous set.
\end{prop}

\begin{proof}  Let $A$ be a $\Delta^0_2$ set such that neither $A$ nor
$\overline{A}$ has an infinite $\Pi^0_1$ subset.  (Such a set is easy to
construct using a wait--and--see argument with $\zerop$ oracle.)  By a result  
of Jockusch (Proposition~2.1 of \cite{Jo1}), there exists a computable stable  
$2$--coloring $\color{C}_A$ of $\ktup{\omega}{2}$ such that every infinite set  
$B$ which is homogeneous for $\color{C}_A$ has the property that either $B
\subseteq A$ or $B \subseteq \overline{A}$.  This result then follows
immediately.
\end{proof}

For our purposes, the notion of a $\Sigma^0_2$ $2$--coloring of
$\ktup{\omega}{n}$  induced by a c.e.~stable $2$--coloring of
$\ktup{\omega}{n+1}$ will be most useful.

\begin{prop}\label{prop:induced}  Given a c.e.~stable $2$--coloring
$\color{C}$ of $\ktup{\omega}{n+1}$, there exists a $\Sigma^0_2$ $2$--coloring  
$\color{P}$ of $\ktup{\omega}{n}$ such that

\begin{enumerate}
\item every infinite set which is homogeneous for $\color{C}$ is also
homogeneous for $\color{P}$, and
\item if $A$ is an infinite homogeneous set for $\color{P}$, then there
exists an infinite set $B \subseteq A$ such that $B$ is homogeneous for
$\color{C}$ and $B \T A \join \zerop$.
\end{enumerate}
\end{prop}

\begin{proof}  Let $\color{C}$ be a c.e.~stable $2$--coloring of
$\ktup{\omega}{n+1}$ which is given as a red--blue coloring of pairs, where  
the set of red pairs of $\color{C}$ is  c.e.~and the set of blue pairs of
$\color{C}$ is co--c.e.  We induce a $\Sigma^0_2$ $2$--coloring $\color{P}$ of  
$\ktup{\omega}{n}$ as follows.  If $a_1 < \dots < a_n$, then let $\{a_1,
\dots, a_n\}$ be $\color{P}$--blue if and only if there exists $s_0$ such that  
for all $s \ge s_0$, $\{a_1, \dots, a_n, s\}$ is $\color{C}$--blue.  By
stability, it follows that $\{a_1, \dots, a_n\}$ is $\color{P}$--red if and  
only if there exists $s_0$ such that for all $s \ge s_0$, $\{a_1, \dots, a_n,  
s\}$ is $\color{C}$--red.

Let $A$ be an infinite homogeneous set for $\color{C}$, and let $a_1, \dots,  
a_n \in A$ with $a_1 < \dots < a_n$.  If $A$ is $\color{C}$--red homogeneous,  
then for all $a_{n+1} > a_n$ with $a_{n+1} \in A$, we have $\{a_1, \dots,
a_n, a_{n+1}\}$ $\color{C}$--red.  Since $A$ is infinite, $\{a_1, \dots,
a_n\}$ is $\color{P}$--red by definition, and so $A$ is homogeneous for
$\color{P}$.  The argument is analogous if $A$ is $\color{C}$--blue
homogeneous.

Next, let $A$ be an infinite homogeneous set for $\color{P}$.
It is easy to show that  there exists an infinite $B
\subseteq A$ with $B \T A \join \zerop$ such that
 $B$ is homogeneous for $\color{C}$, as in the proof of
\propref{prop:stabled2}.
\end{proof}

\begin{cor}\label{cor:stabledel4}  Every c.e.~stable $2$--coloring of
$\ktup{\omega}{3}$ has an infinite $\Delta^0_4$ homogeneous set.
\end{cor}

\begin{proof}  Let $\color{C}$ be a c.e.~stable $2$--coloring of
$\ktup{\omega}{3}$, and let $\color{P}$ be the induced $\Sigma^0_2$
$2$--coloring of $\ktup{\omega}{2}$ from \propref{prop:induced}.  Relativizing  
\thmref{thm:pi2re} to $\zerop$, we see that the $\Sigma^0_2$ $2$--coloring
$\color{P}$ has an infinite $\Pi^0_3$ homogeneous set $A$.  By
\propref{prop:induced}, there is an infinite $\color{C}$--homogeneous set $B  
\subseteq A$ with $B \T A \join \zerop$, and hence $B$ is $\Delta^0_4$.
\end{proof}

\begin{cor}\label{cor:redel4}  Every c.e.~$2$--coloring of $\ktup{\omega}{3}$  
has an infinite $\Delta^0_4$ homogeneous set.
\end{cor}

\begin{proof}  Let $A$ be a $\low{2}$ cohesive set (i.e.,
$\jump{A}{\prime\prime} \T \zeropp$); the existence of such a set is
guaranteed by a result of Jockusch and Stephan \cite{JoSt1} (Theorem~2.5),
whose proof is corrected in \cite{JoSt2}.  If we work
inside $A$, then all c.e.~$2$--colorings of $\ktup{\omega}{3}$ appear
stable.  More
precisely, we have the following.

\begin{lemma}\label{lem:cohstable}  If $\color{C}$ is a c.e.~$2$--coloring
of $\ktup{\omega}{3}$ and $A$ is a cohesive set, then $\widehat{\color{C}} =  
\color{C}\restriction \ktup{A}{3}$ is c.e.~in $A$ and stable.
\end{lemma}

\begin{proof}  The first statement is obvious.  Let $\color{C}$ be given as  
a red--blue coloring of $\ktup{\omega}{3}$ such that the set of red triples of  
$\color{C}$ is c.e.~and the set of blue triples of $\color{C}$ is co--c.e.   
To see that $\widehat{\color{C}}$ is stable, let $a, b \in A$ with $a<b$.
Then $W = \{c \mid \{a,b,c\} \mbox{ is red}\}$ is c.e.  Since $A$ is
cohesive, $A \subseteq^* W$ or $A \subseteq^* \overline{W}$, which proves
that $\widehat{\color{C}}$ is stable.
\end{proof}

Let $\color{C}$ be a c.e.~$2$--coloring of $\ktup{\omega}{3}$, and let
$\widehat{\color{C}} = \color{C} \restriction \ktup{A}{3}$, so that
$\widehat{\color{C}}$ is c.e.~in $A$ and stable.  By \corref{cor:stabledel4}  
relativized to $A$, there exists an infinite set $B$ which is $\Delta^0_4$
relative to $A$ and homogeneous for $\color{C}$.  Since $A$ is $\low{2}$, $B$  
is $\Delta^0_4$ as well.
\end{proof}

Note that this corollary increases our knowledge of c.e.~$2$--colorings of
$\ktup{\omega}{3}$, as before we only knew, by relativizing
\thmref{thm:basicjockusch}, that such $2$--colorings have a $\Pi^0_4$ infinite  
homogeneous set.  We generalize this result to all c.e.~$2$--colorings of
$\ktup{\omega}{n}$, for $n \ge 3$, in \secref{sec:general}.

Finally, we note that just as every c.e.~stable $2$--coloring of
$\ktup{\omega}{n+1}$ induces a $\Sigma^0_2$ $2$--coloring of
$\ktup{\omega}{n}$, every $\Sigma^0_2$ $2$--coloring of $\ktup{\omega}{n}$
induces a c.e.~stable $2$--coloring of $\ktup{\omega}{n+1}$.

\begin{prop}  Given a $\Sigma^0_2$ $2$--coloring $\color{P}$ of
$\ktup{\omega}{n}$, there exists a c.e.~stable $2$--coloring $\color{C}$ of  
$\ktup{\omega}{n+1}$ such that

\begin{enumerate}
\item every infinite set which is homogeneous for $\color{C}$ is also
homogeneous for $\color{P}$, and
\item if $A$ is an infinite homogeneous set for $\color{P}$, then there
exists an infinite set $B \subseteq A$ such $B$ is homogeneous for
$\color{C}$ and $B \T A \join \zerop$.
\end{enumerate}
\end{prop}

\begin{proof}  Let $\color{P}$ be a $\Sigma^0_2$ $2$--coloring of
$\ktup{\omega}{n}$, given as a subset of $\ktup{\omega}{n}$.  Since
$\color{P} \m \fin{\null} = \{e \mid W_e \mbox{ is finite}\}$, fix a
computable function $f$ such that $\{a_1, \dots, a_n\} \in \color{P}$ if and  
only if $W_{f(a_1, \dots, a_n)}$ is finite.  Define a $2$--coloring
$\color{C}$ of $\ktup{\omega}{n+1}$ as a subset of $\ktup{\omega}{n+1}$, as  
follows.  If $a_1< \dots < a_n < a_{n+1}$, then
\[\{a_1, \dots, a_n, a_{n+1}\} \in \color{C} \iff (\exists d \ge
a_{n+1})[W_{f(a_1, \dots, a_n),d} \ne W_{f(a_1, \dots, a_n),d+1}].\]

The $2$--coloring $\color{C}$ is clearly c.e.  To see that $\color{C}$ is
stable, fix $a_1, \dots, a_n$ with $a_1 < \dots < a_n$.  If $\{a_1, \dots,
a_n\} \in \color{P}$, then $W_{f(a_1, \dots, a_n)}$ is finite, and hence for  
sufficiently large numbers $s$, $W_{f(a_1, \dots, a_n),s} = W_{f(a_1, \dots,  
a_n),s+1}$.  Thus $\{a_1, \dots, a_n, a_{n+1}\} \in \overline{\color{C}}$ for  
$a_{n+1}$ sufficiently large.  If $\{a_1, \dots, a_n\} \in
\overline{\color{P}}$, then $W_{f(a_1, \dots, a_n)}$ is infinite, and hence  
$W_{f(a_1, \dots, a_n),s} \ne W_{f(a_1, \dots, a_n),s+1}$ for infinitely many  
numbers $s$.  Hence $\{a_1, \dots, a_n, a_{n+1}\} \in \color{C}$ for all
$a_{n+1} > a_n$.

Let $A$ be an infinite homogeneous set for $\color{C}$.  We show that $A$ is  
homogeneous for $\color{P}$.  Let $a_1, \dots, a_n \in A$ with $a_1 < \dots  
< a_n$.  If $\ktup{A}{n+1} \subseteq \color{C}$, then for all $c \in A$ with  
$c > a_n$, $\{a_1, \dots, a_n, c\} \in \color{C}$.  Since $A$ is infinite,
the definition of $\color{C}$ implies that $W_{f(a_1,\dots, a_n)}$ is
infinite, and hence $\{a_1, \dots, a_n\} \in \overline{\color{P}}$.  Thus
$\ktup{A}{n} \subseteq \overline{\color{P}}$.  If $\ktup{A}{n+1} \subseteq
\overline{\color{C}}$, then for all $c \in A$ with $c > a_n$, $\{a_1, \dots,  
a_n,c\} \in \overline{\color{C}}$.  Since $A$ is infinite, the definition of  
$\color{C}$ implies that $W_{f(a_1, \dots, a_n)}$ is finite, and hence
$\{a_1, \dots a_n\} \in \color{P}$.  Thus $\ktup{A}{n} \subseteq \color{P}$.

Let $A$ be an infinite homogeneous set for $\color{P}$.  We show that there  
exists an infinite set $B \subseteq A$ with $B \T A \join \zerop$ and $B$
homogeneous for $\color{C}$.  First assume that $\ktup{A}{n} \subseteq
\color{P}$.  Let $b_1, \dots, b_{n}$ be the least $n$ elements of $A$.  Since  
$\{b_1, \dots, b_{n}\} \in \color{P}$, $W_{f(b_1, \dots, b_{n})}$ is finite,  
and hence for sufficiently large numbers $c$, $\{b_1, \dots, b_{n}, c\} \in  
\overline{\color{C}}$.  Let $b_{n+1} = (\mu x \in A)[x > b_{n} \land \{b_1,  
\dots, b_{n},x\} \in \overline{\color{C}}]$.  Then for all $1 \le i_1 < \dots  
< i_{n} \le n+1$, $\{b_{i_1}, \dots, b_{i_n}\} \in \color{P}$ and hence, as  
before, we can let $b_{n+2}$ be the least $x \in A$ such that $x > b_{n+1}$  
and for all $1 \le i_1 < \dots < i_{n} \le n+1$, $\{b_{i_1}, \dots, b_{i_n},  
x\} \in \overline{\color{C}}$.  We can continue in this fashion to enumerate  
a set $B = \{b_1, b_2, \dots\}$ with $\ktup{B}{n+1} \subseteq
\overline{\color{C}}$.  We have $B \T A \join \zerop$, since $B$ is
enumerated in increasing order and $\overline{\color{C}}$ is $\Pi^0_1$.

Next assume that $\ktup{A}{n} \subseteq \overline{\color{P}}$.  Let $a_1,
\dots, a_{n+1} \in A$ with $a_1 < \dots < a_{n+1}$.  Then $\{a_1, \dots,
a_n\} \in \overline{\color{P}}$, and hence $W_{f(a_1, \dots, a_n)}$ is
infinite.  Thus $\{a_1, \dots, a_n, a_{n+1}\} \in \color{C}$, and hence
$\ktup{A}{n+1} \subseteq \color{C}$.  Thus we take $B = A$.
\end{proof}

\section{$2$--cohesive and $2$--r--cohesive sets}

\subsection{$2$--cohesive sets}  We begin with a study of $2$--cohesive sets.   
We have already noted that the $1$--cohesive sets are exactly the cohesive
sets; similarly, the $1$--r--cohesive sets are exactly the r--cohesive sets.  It  
is easy to prove that for all $n \ge 1$, $n$--cohesive sets and
$n$--r--cohesive sets exist.

\begin{prop}  If $\color{C}_i$ is a $k_i$--coloring of $\ktup{\omega}{n_i}$,  
for $i \in \omega$, then there is an infinite set $D$ such that for all $i$,  
$D$ is almost homogeneous for $\color{C}_i$.
\end{prop}

\begin{proof}    We construct infinite sets $A_0 \supseteq A_1 \supseteq
\dots$ such that for each $i$, $A_i$ is homogeneous for $\color{C}_i$.

Let $A_0$ be an infinite homogeneous set for $\color{C}_0$.  Given $A_0,
\dots, A_i$, let $\widehat{\color{C}_{i+1}} = \color{C}_{i+1} \restriction
\ktup{A_i}{n_i}$, and let $A_{i+1}$ be an infinite homogeneous set for
$\widehat{\color{C}_{i+1}}$, so that $A_i \supseteq A_{i+1}$.  We then form a  
set $D$ which is a diagonal intersection of the $A_i$'s.  Let $d_0 \in A_0$,  
and given $d_0, \dots, d_i$, let $d_{i+1} \in A_{i+1} - \{d_0, \dots,
d_i\}$.  Then $D = \{d_0, d_1, \dots\}$ is almost homogeneous for each
$\color{C}_i$, since for all $j \ge i$, $d_j \in A_i$.
\end{proof}

It follows immediately that $n$--cohesive and $n$--r--cohesive sets exist.   
Such a construction of an $n$--cohesive or $n$--r--cohesive set
is quite nonconstructive,  as the set $D$ constructed above is not
obviously arithmetical.  It is possible to show, however, that
for $n \ge 1$, an arithmetical $n$--cohesive,  and hence an
arithmetical $n$--r--cohesive, set exists.

The existence of a maximal set, a c.e.~set whose complement is cohesive,
shows that a $\Pi_1^0$ $1$--cohesive set exists.  In
the other direction, it is clear from the fact that no c.e.~set
is r--cohesive and \thmref{thm:basicjockusch} that no $n$--r--cohesive
set is $\Sigma^0_n$ for any $n \geq 1$.  Thus the next result is
the best possible result in terms of the arithmetical hierarchy
for existence of $2$--cohesive sets.

\begin{thm}\label{thm:pi22co}  There exists a $\Pi_2^0$ $2$--cohesive set $A$.
\end{thm}

\begin{proof}  Under a suitable coding of pairs of natural numbers, we can
enumerate the c.e.~$2$--colorings $W_0, W_1, \dots$ of pairs, where each $W_i$  
is a c.e.~subset of $\ktup{\omega}{2}$.  Let each $2$--coloring $W_i$ be a
red--blue coloring of $\ktup{\omega}{2}$, where $\{x,y\}$ ($x \ne y$) is {\em  
red\/} if $\{x,y\} \in W_i$, and $\{x,y\}$ is {\em blue\/} otherwise.  If
$\{x,y\} \in W_i$, we shall say that the $i$--color of $\{x,y\}$ is red; if  
$\{x,y\}\notin W_i$, we shall say that the $i$--color of $\{x,y\}$ is blue.

We first recall the maximal set construction, which is just a construction
of a $\Pi_1^0$ ($1$--)cohesive set.  The $\Pi_1^0$ cohesive set is constructed  
by a movable marker construction, and each number is labelled with an
$e$--state at stage $s$.  Given a number $x$, the {\em $e$--state of $x$ at  
stage $s$\/} is $\sigma(e,x,s) = \{i \mid i \le e \land x \in W_{i,s}\}$.
Each $e$--state is identified with an ($e+1$)--digit binary number, such that  
the $i$th bit ($i \le e$, read from left to right) is $0$ if $x \notin W_i$  
and $1$ if $x \in W_i$.  The $e$--states are ordered lexicographically.  In  
the construction, the $e$th marker moves in order to maximize its $e$--state.   
Since for each $x$, there are only finitely many $e$--states, each marker can  
move only finitely often, and it will follow that the set $A$ which is
defined by the marker construction must have the property that for every $e$,  
either $A \subseteq^* W_e$ or $A \subseteq^* \overline{W_e}$; i.e., for
every $e$ there will be a finite set $F$ and a fixed $e$--state $\sigma$ such  
that for all $x \in A-F$, $x$ has $e$--state $\sigma$.

We want to use the idea of the maximal set construction, as well as the idea  
of the construction of an infinite $\Pi_2^0$ homogeneous set for a given
c.e.~$2$--coloring of $\ktup{\omega}{2}$ (\thmref{thm:pi2re}), modified to
consider all c.e.~$2$--colorings of $\ktup{\omega}{2}$.  The construction will  
be a movable marker construction with a $\zerop$ oracle and will result in
an increasing sequence $\{a_n\}_{n\in \omega}$ of numbers.  During the
construction, as in the proof of \thmref{thm:pi2re}, each $a_e$ will be
colored an $i$--color with respect to the $2$--coloring $W_i$, for all $i \le  
e$.  We want the sequence $\{a_n\}_{n \in \omega}$ to have the property that,  
for all $i \le e$ and all $m > e$, the $i$--color of $a_e$ is the same as the  
$i$--color of $\{a_e,a_m\}$.  We will use $e$--states to keep track of the
$i$--colors of $a_e$ and, as in the maximal set construction, maximize the
$e$--state of $a_e$ in order to ensure that, from some point on, all $a_n$'s  
have the same $e$--state.

We denote the position of marker $\Lambda_e$ at the beginning of stage $s$
by $a_e^s$.    Given $a_e^s$, denote the $e$--state of $a_e^s$ by an
($e+1$)--digit binary number, such that the $i$th bit ($i\le e$, read from
left to right) of the $e$--state is $0$ (respectively $1$) if the $i$--color of  
$a_e^s$ is red (respectively blue).  For a fixed $e$, the $i$--state of
$a_e^s$, where $i<e$, is defined in the obvious way.  As in the maximal set  
construction, $e$--states are read from left to right and ordered
lexicographically.  As an example, if the $2$--state of $x$ at $s$ is $100$,  
then the $0$--color of $x$ at $s$ is blue, and the $1$--color and $2$--color of  
$x$ at $s$ are both red.  Let $\sigma(x,e,s)$ denote the $e$--state of a
marked number $x$ at stage $s$.  Let $\sigma$ and $\tau$ be $e$--states, with  
$\tau > \sigma$, and $\sigma$ the $e$--state of some number with marker
$\Lambda_i$, and $\tau$ the $e$--state of some number with marker $\Lambda_j$,  
$e \le i \le j$.  Let $n_0$ be the least $n$, $0 \le n \le e$, such that
$\sigma$ and $\tau$ differ in the $n$th bit.  Then we say that $\Lambda_i$
wants to improve its $e$--state for the sake of $W_{n_0}$.

In this construction, a marker can move for two reasons:  to improve its
$e$--state as in the maximal set construction, or as in the proof of
\thmref{thm:pi2re}.  As before, we initially assume the correct $j$--color, $j  
\le i$, for $a_i$ is red.  In the proof of \thmref{thm:pi2re}, a number $c$  
was called $k$--acceptable at $s$ if for all $i < k$, $a_i^s$ is defined,
$a_i^s < c$, and $\{a_i^s,c\}$ is red (in the given fixed $2$--coloring) if  
$a_i^s$ is red at stage $s$.  Thus $c$ is acceptable (with respect to the
given $2$--coloring) at $s$ to the entire initial segment $a_0^s, \dots,
a_{k-1}^s$.  At stage $s$ of the construction, there are exactly $n(s)$
numbers $a_0^s,\dots a_{n(s)-1}^s$ assigned to markers $\Lambda_0, \dots,
\Lambda_{n(s)-1}$.  The largest $j \le n(s)$ is found such that there exists  
a number $c$ which is $j$--acceptable.  If $j< n(s)$, then $c$ is acceptable  
to the entire initial segment $a_0^s, \dots, a_{j-1}^s$, and the color of
$a_j^s$ is changed from red to blue, which makes $c$ acceptable to $a_j^s$ as  
well.  If $j=n(s)$, then $c$ is acceptable to all numbers with a marker and  
the marker $\Lambda_{n(s)}$ is assigned.

We need to prioritize our $2$--colorings $W_0 > W_1 > W_2 > \dots$ (with $W_i  
> W_j$ meaning that $W_i$ has higher priority than $W_j$).  We need a notion  
of acceptability which respects the priority ranking of the $W_i$'s and the  
fact that a number assigned to a marker $\Lambda_e$ has an $e$--state of
colors assigned to it.  We define a notion of $(k_1,k_2)$--acceptability as  
follows.

\begin{definition}  Let $k_2 \le k_1$.  A number $c$ is {\em
$(k_1,k_2)$--acceptable at stage $s$\/} if
\begin{enumerate}
\item for all $i<k_1$, $a_i^s$ is defined, $a_i^s<c$, and for all $j \le i$,  
if $a_i^s$ is $j$--red, then $\{a_i^s,c\}$ is $j$--red, and
\item if $k_2 \ne 0$, then $a_{k_1}^s$ is defined, $a_{k_1}^s<c$, and for
all $j < k_2$, if $a_{k_1}^s$ is $j$--red, then $\{a_{k_1}^s,c\}$ is $j$--red.
\end{enumerate}
\end{definition}
To clarify this notion, note that a number $c$ is $(k_1,k_2)$--acceptable at  
$s$ if $c$ is acceptable to the entire initial segment $a_0^s, \dots,
a_{k_1-1}^s$ (i.e., for all $i<k_1$, $c$ is acceptable to the entire
$i$--state of $a_i^s$), and if $a_{k_1}^s$ is defined, then $c$ is acceptable  
to the $j$--colors of $a_{k_1}^s$, for $j<k_2$, but not acceptable to the
$k_2$--color of $a_{k_1}^s$.

Finally, say that a number $c$ is {\em free at $s$\/} if, prior to stage
$s$, it has not been the position of any marker, and $c \ge s$.

The idea of this construction will be as follows.  Suppose at stage $s$
there are exactly $n(s)$ numbers $a_0^s, \dots, a_{n(s)-1}^s$ assigned to
markers $\Lambda_0, \dots, \Lambda_{n(s)-1}$.  We want to find a number which  
is acceptable to all colors of as long an initial segment of $a_0^s, \dots,  
a_{n(s)-1}^s$ as possible, and acceptable to as many of the colors of the
next $a_j^s$ (if it exists) as possible, respecting the priority of the
$W_i$'s.

\noindent {\em Construction.}

\noindent {\em Stage $2s$.}  Assume inductively that there is a number
$n(2s)$ such that the markers having a position are exactly the $\Lambda_i$,  
$i<n(2s)$.

\setcounter{case}{0}
\begin{case}  There exists a number $c$ which is free and
$(n(2s),0)$--acceptable at $2s$.
\end{case}

Attach $\Lambda_{n(2s)}$ to the least such number $c$, and color all
$n(2s)+1$ colors of $c$ red (i.e., the $n(2s)$--state of $c$ is
$\underbrace{00\dots 0}_{n(2s)+1}$).

\begin{case}  Otherwise.
\end{case}

Let $j(2s)$ be the largest number $j$ for which there exists a $k \le j$ and  
some number which is free and $(j,k)$--acceptable at $2s$.  Given $j(2s)$,
let $k(2s)$ be the greatest such $k$.  Such numbers exist because every
number is $(0,0)$--acceptable at $2s$.  Note that $j(2s) < n(2s)$.  We
\begin{enumerate}
\item change the $k(2s)$--color of $a_{j(2s)}^{2s}$ (necessarily
from red to blue),
\item let the $i$--color of
$a_{j(2s)}^{2s}$ be red, for $k(2s) < i \le  j(2s)$, and
\item detach marker $\Lambda_i$, for $j(2s) < i < n(2s)$.
\end{enumerate}

\noindent {\em Stage $2s+1$.}  Let $n(2s+1)$ be such that the markers having  
a position are exactly the $\Lambda_i$, $i < n(2s+1)$.  Choose the least
$i<n(2s+1)$ such that for some $j$, $i<j<n(2s+1)$ and the $i$--state of
$a_j^{2s+1}$ is greater than the $i$--state of $a_i^{2s+1}$.  (If $i$ fails to  
exist, then proceed to the next stage.)  For this $i$, choose the least such  
$j$.  Move marker $\Lambda_i$ to $a_j^{2s+1}$, and move as many markers
$\Lambda_k$, $k >i$, as possible, preserving their order, to elements
$a_{\hat{k}}^{2s+1}$, $\hat{k} > j$.  (There will not be enough marked
positions; move only those possible and detach the rest of the markers.)

As always, at each stage, any unmentioned markers are left unchanged.  Note  
that only a $\zerop$ oracle is required for the construction.

\begin{lemma}  For any $n$, $\lim_s a_n^s \downarrow$, and $\lim_s
\sigma(a_n,n,s)$ exists.
\end{lemma}

\begin{proof}  Assume inductively that the lemma holds for all $k<n$ and
prove it for $n$.  Let $s_0$ be the least even stage such that for all
$k<n$, $\Lambda_k$ has position $a_k$ at all stages $s_1 \ge s_0$ and
$\sigma(a_k,k,s_1) = \lim_s \sigma(a_k,k,s) = \sigma(a_k,k)$ for all $s_1 \ge  
s_0$.  Then either $\Lambda_n$ is already attached to some number $c =
a_n^{s_0}$ at the beginning of stage $s_0$, or $\Lambda_n$ becomes attached  
to a number $c$ through Case~1 in this even stage of the construction.  After  
stage $s_0$, note that whenever the $n$--state of a number associated with
$\Lambda_n$ changes, the $n$--state actually strictly  increases.

To see this, let $s \ge s_0$ be an even stage at which the $n$--state of
$a_n^s$ changes.  At stage $s$, there exist free numbers which are
$(n,0)$--acceptable, since otherwise, for some $k<n$, either $\Lambda_k$
becomes detached, or the $k$--state associated with $\Lambda_k$ changes.  So,  
it must be the case that for some $k$, $0 \le k < n$, there exist
$(n,k)$--acceptable numbers but not $(n,k+1)$--acceptable numbers.  It follows  
by the definition of acceptability that the $k$--color of $a_n^s$ is red, but  
there do not exist any free $(n,k)$--acceptable numbers making a $k$--red pair  
with $a_n^s$.  So, all free $(n,k)$--acceptable numbers must make a $k$--blue  
pair with $a_n^s$.  Then the $k$--color of $a_n^s$ is changed from red to
blue, which causes the $n$--state of $a_n^s$ to strictly increase.  Also,
whenever marker $\Lambda_n$ moves at an odd stages after $s_0$, the $n$--state  
strictly increases, by construction.

Since there are only finitely many $n$--states, there exists $s_1 \ge s_0$
such that for all $s \ge s_1$, $a_n^s = a_n^{s_1} = a_n$.  After stage $s_1$,  
if the $n$--state of $a_n$ changes because of an even stage in the
construction, then again, the $n$--state must strictly increase.  So, $\lim_s  
\sigma(a_n,n,s)$ exists, since there are only finitely many $n$--states.
\end{proof}

\begin{lemma}\label{lem:stabcolor}  For each $i$, $\lim_e (i\mbox{--color of  
}a_e)$ exists.
\end{lemma}

\begin{proof}  Fix $i$ and assume the lemma for all $j<i$.  Choose $e_0$
such that for all $e \ge e_0$, for all $j<i$, the $j$--color of $a_e$ is the  
same as the $j$--color of $a_{e_0}$.  Assume the lemma does not hold for $i$.   
Then we can choose $e_0<e_1<e_2$ such that the $i$--color of $a_{e_1}$ is red  
and the $i$--color of $a_{e_2}$ is blue.  Choose $s$ such that for all $k \le  
e_2$, $a_k^s = a_k$ and $\sigma(a_k,k,s) = \sigma(a_k,k)$.  But then
$\sigma(a_{e_1},i) < \sigma(a_{e_2},i)$, and so some marker $\Lambda_r$, $r  
\le e_1$, moves at stage $s+1$, which is a contradiction.
\end{proof}

\begin{lemma}\label{lem:stabpair}  If $i<j$, then for all $k \le i$, the
pair $\{a_i,a_j\}$ has the same $k$--color as the eventual $k$--color of $a_i$.
\end{lemma}

\begin{proof}  Let $i<j$ and fix $k \le i$.  When $a_j$ first is the
position of $\Lambda_j$, $a_j$ must be $(j,0)$--acceptable.  The $k$--color of  
$a_i$ cannot later be changed at even stages, since otherwise $a_j$ loses
its marker and never regains it.  The $k$--color of $a_i$ cannot be changed at  
an odd stage, because $a_j$ would no longer be associated with $\Lambda_j$.   
By definition of $a_j$ being $(j,0)$--acceptable, we see that if the
$k$--color of $a_i$ is red, then the $k$--color of $\{a_i,a_j\}$ is red.  If  
the $k$--color of $a_i$ is blue, then at the stage $s$ where the $k$--color of  
$a_i = a_i^s$ was made blue, there were no more free $(i,k)$--acceptable
numbers making a $k$--red pair with $a_i$.  However, since there are
infinitely many $(i,k)$--acceptable numbers, they must all make a $k$--blue  
pair with $a_i$.  So, any new number associated with $\Lambda_j$ must make a  
$k$--blue pair with $a_i$.  Since $a_j$ is such a number, $\{a_i,a_j\}$ is
blue.
\end{proof}

Let $A = \{a_i \mid i \in \omega\}$.  Note that $A$ is infinite because
$\{a_i\}_{i \in \omega}$ is an increasing sequence.  Also note that $A$ is
$\Pi_2^0$ since any number which loses a marker may never be a marker
position again (although a number may change from one marker to another).
For this reason, for all $x$,
\[x \in \overline{A} \iff (\exists s)[x \le s \land x \mbox{ has no marker
at } s].\]
Since the construction is computable in $\zerop$, we see that $\overline{A}$  
is $\Sigma_2^0$.  By \lemref{lem:stabcolor} and \lemref{lem:stabpair}, we see  
that for any c.e.~$2$--coloring $W$ of $\ktup{\omega}{2}$, $A$ is almost
homogeneous for $W$.
\end{proof}

We are thus led to the ask the following question.

\begin{ques}  For which $n \ge 1$ do
$\Pi_n^0$ $n$--cohesive sets exist?
\end{ques}
Recall that it is still unknown, even for $n = 3$,
whether every c.e.~$2$--coloring of  $\ktup{\omega}{n}$ has an
infinite $\Pi^0_n$ homogeneous set.  However, in
\secref{sec:general}, we construct a $\Delta^0_{n+1}$ $n$--cohesive set, for  
 $n \ge 2$.

Note that there exists a c.e.~$2$--coloring $\color{C}$ of
$\ktup{\omega}{2}$ such that for all infinite homogeneous sets $A$, $\zerop  
\T A$.  Such a $2$--coloring $\color{C}$ is defined as a subset of
$\ktup{\omega}{2}$ as follows.  If $a < b$, then
\[\{a,b\} \in \color{C} \iff (\exists c \le a)[c \in K \land c \notin K_b].\]
(Here, $\{K_s\}_{s \in \omega}$ is a fixed computable enumeration of the
complete c.e.~set $K$.)  It is easy to see that $\color{C}$ is c.e.~and that  
every infinite homogenous set $A$ for $\color{C}$ is such that $\ktup{A}{2}  
\subseteq \overline{\color{C}}$.  It follows immediately that if $A$ is an
infinite homogeneous set, then $K\T A$.

The following proposition is immediate.

\begin{prop}\label{prop:2coabovezp} If $A$ is a $2$--cohesive set, then
$\zerop \T A$.
\end{prop}
We next apply \corref{cor:nopi2andjump} to $\Pi_2^0$ $2$--cohesive sets.

\begin{thm}\label{thm:pi22cozpp}  Every $\Pi_2^0$ $2$--cohesive set is of
degree $\td{\jump{0}{\prime \prime}}$.
\end{thm}

\begin{proof}  Let $M$ be a $\Pi_2^0$ $2$--cohesive set.
\corref{cor:nopi2andjump} says that there exists a computable $2$--coloring of  
$\ktup{\omega}{2}$ such that for all infinite homogeneous sets $A$, if $A$ is  
$\Pi_2^0$, then $\zeropp \T A \join \zerop$.  Since $M$ differs finitely
from a $\Pi_2^0$ set $A$ homogeneous for this partition, we have $\zeropp \T  
M \join \zerop$, and hence $\zeropp \T M$, since $\zerop \T M$ by
\propref{prop:2coabovezp}.  Thus $\zeropp \eT M$, since $M$ is $\Pi_2^0$.
\end{proof}

It is natural to ask whether every $2$--cohesive set has degree at least
$\td{\jump{0}{\prime\prime}}$.  To answer this question requires the
following impressive theorem of Seetapun.

\begin{thm}[Seetapun, \cite{SeSl} (Theorem~2.1)]\label{thm:seet}  For any
computable
$2$--coloring of $\ktup{\omega}{2}$ and any noncomputable sets $C_0, C_1,
\dots$, there exists an infinite homogeneous set $A$ such that for all $n$,  
$C_n \not\T A$.
\end{thm}
Seetapun's theorem implies that for every computable $2$--coloring of
$\ktup{\omega}{2}$, there exists an infinite homogeneous set $A$ with $\zerop  
\not\T A$.  We will require the following relativized version of Seetapun's  
theorem, namely, given a $2$--coloring $\color{C}$ of $\ktup{\omega}{2}$ with  
$\color{C} \T B$, and a set $D \not\T B$, there exists an infinite
homogeneous set $A$ with $D \not\T A \oplus B$.

\begin{thm}\label{thm:kps} The following are equivalent for every degree $\td{d}$.

\begin{enumerate}
\item Every $2$--cohesive set has degree at least $\td{d}$.
\item $\td{d} \le \td{\jump{0}{\prime}}$.
\end{enumerate}
\end{thm}

\begin{proof}  ((2) $\implies$ (1))  This follows immediately from
\propref{prop:2coabovezp}.

((1) $\implies$ (2))  We prove the contrapositive.  Let $D$ be a set with $D  
\not\T \zerop$.  We use \thmref{thm:seet} to construct a $2$--cohesive set
$A$ such that $D \not\T A$.

Let $\color{C}_0, \color{C}_1, \dots$ be a list of all $\zerop$--computable  
$2$--colorings of $\ktup{\omega}{2}$; this list contains all
c.e.~$2$--colorings of $\ktup{\omega}{2}$.  We construct infinite sets $A_0  
\supseteq A_1 \supseteq \dots$ such that for each $i$, $A_i$ is homogeneous  
for $\color{C}_i$, and $D \not\T A_0 \join \cdots \join A_i \join \zerop$.

First apply the relativized version of \thmref{thm:seet} for the
non--$\zerop$--computable set $D$ to $\color{C}_0$ to obtain an infinite
homogeneous set $A_0$ with $D \not\T A_0 \join \zerop$.  Next, assume $A_0,  
\dots, A_n$ are given such that for $0 \le i \le n$, $A_i$ is an infinite
homogeneous set for $\color{C}_i$,  $D \not\T A_0 \join \cdots \join A_n
\join \zerop$, and $A_0 \supseteq A_1 \supseteq \dots \supseteq A_n$.  Let
$\widehat{\color{C}_{n+1}} = \color{C}_{n+1}\restriction \ktup{A_n}{2}$, so  
that $\widehat{\color{C}_{n+1}} \T A_n \join \zerop \T A_0 \join \cdots \join  
A_n \join \zerop$.  Apply the relativized version of \thmref{thm:seet} for
$D$ to $\widehat{\color{C}_{n+1}}$ to obtain an infinite homogeneous set
$A_{n+1}$ for $\color{C}_{n+1}$ such that $D \not\T A_0 \join \cdots \join
A_n \join A_{n+1} \join \zerop$.  Note that the proof of \thmref{thm:seet}
can be modified to obtain $A_{n+1} \subseteq A_n$ since $D \not\T A_n$.  We  
thus have for all $n$, $D \not\T A_0 \join \cdots \join A_n \join \zerop$.

By a slight modification of the proof of the Kleene--Post--Spector Theorem on  
exact pairs (see \cite{So1},  Theorem~VI.4.2), there exist sets
$\widehat{A_i}$, $i \in \omega$, such that for all $i$, $\widehat{A_i} =^*
A_i$ and $D \not\T \bigoplus\{\widehat{A_i} \mid i \in \omega\}$.  (The proof  
of the Kleene--Post--Spector Theorem needs to  be slightly modified, as the  
original proof requires $A_0 <_T A_1 <_T \dots$.  However, the hypothesis
that for all $n$, $D \not\T A_0 \join \cdots \join A_n \join \zerop$, is
enough to prove the desired result.)

The technique of ``diagonal intersection'' now produces an
infinite set $C$ such that  $C \subseteq^* \widehat{A_i}$ for all $i$, and $C 
\T \bigoplus\{\widehat{A_n} \mid n \in \omega\}$.  (Let $C =
\{c_n \mid n \in \omega\}$, where $c_n$  is the least element of
$\cap_{i \leq n} \widehat{A_i} - \{c_i \mid i < n \}$.  Note that
$\cap_{i \leq n} \widehat{A_i}$  is infinite because it differs
only finitely from $A_n$.)  Thus,  $C$ is $2$--cohesive and $D
\not\T C$, since $C \T \bigoplus\{\widehat{A_i}  \mid i \in
\omega\}$ and $D \not\T \bigoplus\{\widehat{A_i} \mid i \in
\omega\}$. \end{proof}

It follows that there exists a $2$--cohesive set $A$ with $\zeropp \not\T
A$.  Finally, we give a result regarding the jump of a
$2$--cohesive  set.

\begin{prop}\label{prop:high2co}  If $A$ is $2$--cohesive, then
$\jump{0}{(3)} \T \jump{A}{\prime}$.
\end{prop}

\begin{proof}  Assume that $A$ is $2$--cohesive.  We show that
$p_A$, the  function which enumerates $A$ in increasing order,
dominates all  $\zerop$--computable functions; i.e., for every
$\zerop$--computable function  $f$ and all sufficiently large
$n$, $p_A(n) \ge f(n)$.   Without loss of generality,
we may restrict our attention to such $f$ which are increasing,
and we may assume that $0 \notin A$.  Let $f \T \zerop$,
and let $g$ be a computable function such that for all $x$, $f(x)
= \lim_s g(x,s)$.  Consider  the following c.e.~$2$--coloring
$\color{C}$ of $\ktup{\omega}{2}$, defined  as a c.e.~subset of
$\ktup{\omega}{2}$, \[\color{C} = \{\{a,b\} \mid a < b \land
(\exists s)[b \le g(a,s)]\}.\] Since $A$ is $2$--cohesive, there
is a finite set $F$ such that $A-F$ is homogeneous for
$\color{C}$.  Since $A-F$ is infinite and for all $x$,  $\lim_s
g(x,s)$ exists, it follows that $\ktup{A-F}{2} \subseteq
\overline{\color{C}}$.  It then follows that for all $n$
sufficiently large  and for all $s$, $p_A(n+1) > g(p_A(n),s)$.
Hence $p_A(n+1) > f(p_A(n)) \ge  f(n+1)$, for all sufficiently
large $n$,  as $f$ is increasing and $p_A(n) \ge n+1$ because $0
\notin A$.   It follows that $p_A$ dominates all $\zerop$--computable
functions, and hence by a theorem of Martin \cite{Ma2}
(Lemma~1.1), relativized to $\zerop$, $\jump{0}{(3)} \T \jump{(A
\oplus \zerop)}{\prime}$.  Since $\zerop \T A$ by
\propref{prop:2coabovezp}, it follows that $\jump{0}{(3)} \T
\jump{A}{\prime}$. \end{proof}

\subsection{$2$--r--cohesive sets}

We now turn our attention to $2$--r--cohesive sets, infinite sets which are  
almost homogeneous for every computable $2$--coloring of $\ktup{\omega}{2}$.   
We begin with the following result concerning the jumps of degrees
of $2$--r--cohesive sets,  part of which is the analogue of
\propref{prop:high2co}.
\begin{thm}\label{thm:high2rco} If $A$ is
$2$--r--cohesive, then $\zeropp <_T
\jump{A}{\prime}$.
\end{thm}

\begin{proof}  The result that if $A$ is $2$--r--cohesive then  $\zeropp \T
\jump{A}{\prime}$ is due to Stephan \cite{St}.  We
prove this part first.  By  \cite{Ma2}  (Lemma~1.1), it
suffices to show that $p_A$ dominates all computable functions
$f$.  This is done by a proof that is parallel to that of
\propref{prop:high2co} but considering for each increasing
computable function $f$ the computable $2$--coloring of pairs
$\color{C}_f$ in which the pair $\{x,y\}$ with $x < y$ is colored
red if and  only if $y \leq f(x)$.

We now show that if $A$ is $2$--r--cohesive, then $\jump{A}{\prime}
\not\T \zeropp$.  Assume for a contradiction that $A$ is
$2$--r--cohesive and $\jump{A}{\prime} \T \zeropp$.  We define a
computable red--blue coloring $\color{C} = \varphi_a$ of $\ktup{\omega}{2}$  
for which $A$ is not almost homogeneous;
i.e., $\color{C}$ should satisfy the requirements \begin{align*}
R_{2n}&: (\exists a, b \in A)[a > b \ge n \land \{a,b\} \mbox{ is
red}]\\ R_{2n+1}&: (\exists a, b \in A)[a > b \ge n \land
\{a,b\} \mbox{ is blue}]. \end{align*}

By the recursion theorem, we may use the index $a$ of
$\color{C}$ in our construction of $\color{C}$.  Note that $\{i \mid
R_i \mbox{ is satisfied}\}$ is a $\Sigma^A_1$, hence
$\zeropp$--computable, set, whose index as such depends effectively on
$a$ (without assuming in advance that $\varphi_a$ is total).
It follows from the iterated Limit Lemma, and the uniformity of its
proof, that we may effectively compute from $a$ an index of a ternary
computable function $f$ such that for all $i$, if $R_i$ is satisfied
then $\lim_t \lim_s f(i,s,t) = 1$, and if $R_i$ is not satisfied then
$\lim_t \lim_s f(i,s,t) = 0$.  Thus the recursion theorem allows us to
use $f$ in our construction of the coloring $\color{C} = \varphi_a$.

The computable $2$--coloring of pairs is defined as follows; to color the
pair $\{s, t\}$, where $s<t$, let $i$ be the least $e \le s$ with $f(e,s,t) =  
0$. (Note that $i$ is an approximation to the least $j \le s$ such
that $R_j$ is not satisfied.)  If $i$ is even, or does
not exist, then color  $\{s,t\}$ red.  Otherwise, color $\{s,t\}$
blue.

If all requirements $R_i$ are satisfied, then $A$ is not $2$--r--cohesive,
and we are done.  So, assume otherwise and let $i$ be the least $e$ such that  
$R_e$ is not satisfied.  Without loss of generality, assume that $i$ is
even, so that for all $a,b \in A$ with $i \le a < b$, $\{a,b\}$ is blue.
Since $R_i$ is not satisfied, $\lim_t \lim_s f(i,s,t) = 0$.  Since all $R_e$,  
$e < i$, are satisfied, we can fix $s_0 > i$ large enough such that when $s  
\ge s_0$, there exists $t_s$ such that for all $t \ge t_s$, if $j<i$ then
$f(j,s,t) = 1$, and $f(i,s,t) = 0$.  Since $A$ is infinite, we can take $s
\ge s_0$ with $s \in A$, and for this $s$, we can take $t > \max\{s, t_s\}$  
with $t \in A$.  Then for all $j < i$, $f(j,s,t) = 1$ and $f(i,s,t) = 0$.  By  
construction $\{s,t\}$ is colored red, which is a contradiction.  Hence
every $2$--r--cohesive set $A$ satisfies $\jump{A}{\prime} \not\T \zeropp$.  
\end{proof}

It is unknown whether \thmref{thm:pi22cozpp} holds for
$2$--r--cohesive sets or whether every $2$--r--cohesive set $A$
satisfies $\jump{0}{(3)} \leq_T \jump{A}{\prime}$. Next, we note that
while \propref{prop:2coabovezp} states that every $2$--cohesive set
has degree at least $\mathbf \zerop$, this result does not hold for
$2$--r--cohesive sets.

\begin{thm}\label{thm:2rcozp}  Let $D$ be noncomputable.  There exists a
$2$--r--cohesive set $A$ such that $D \not\T A$.
\end{thm}

\begin{proof}  The proof is the same as that for the ((1) $\implies$ (2))
direction of \thmref{thm:kps}, with $\zerop$--computable
$2$--colorings of $\ktup{\omega}{2}$ replaced by computable
$2$--colorings of $\ktup{\omega}{2}$, and the $D$ of \thmref{thm:kps}
replaced by a noncomputable $D$.  (Actually, the proof of
\thmref{thm:kps} is a relativization of the proof of this result to
$\zerop$.)
\end{proof}

The following corollary follows immediately from
\thmref{thm:2rcozp}, since  by \propref{prop:2coabovezp},
$\zerop$ is computable in any $2$--cohesive set  $A$.

\begin{cor}\label{cor:diffdegrees}  There exists a
$2$--r--cohesive set $A$ such that no $2$--cohesive set $B$
satisfies $B \leq_T A$.
\end{cor}

The proofs of \thmref{thm:kps} and \thmref{thm:2rcozp} are
highly nonconstructive, in sense that the sets which are constructed
are not obviously arithmetical, even when the given set $D$ is
arithmetical. However, using a result of Cholak, Jockusch, and Slaman,
it is possible to make the set in \thmref{thm:2rcozp} arithmetical
when $D = K$.  The analogous result for single computable
partitions of $\ktup{\omega}{2}$ was obtained by a different method by
Hummel in \cite{Hu2}, Theorem 3.1.

\begin{thm}\label{thm:d32rco} There exists a $\Delta^0_3$
$2$--r--cohesive set $A$ with $\zerop \not\T A$.
\end{thm}

\begin{proof}  It was shown by Cholak, Jockusch, and Slaman
\cite{ChJoSl} that any computable $2$--coloring of
$\ktup{\omega}{2}$ has an infinite $\low{2}$ homogeneous set;
i.e., an infinite homogeneous set $A$ with
$\jump{A}{\prime\prime} \T \zeropp$.  The idea of the proof of
\thmref{thm:d32rco} is to iterate this result from \cite{ChJoSl}
over all computable $2$--colorings of $\ktup{\omega}{2}$ and then
to take an appropriate diagonal intersection. This result is
convenient for iteration when considering infinitely many
computable $2$--colorings of $\ktup{\omega}{2}$ because if $A$ is
$\low{2}$, and $B$ is $\low{2}$ relative to $A$, then $B$ is also
$\low{2}$.  Of course, when this result is applied infinitely
often to construct a set $C$, it is important to know the extent
to which this result holds {\em uniformly\/} in  order to
calculate the complexity of $C$.  The following lemma is a
relativized version of the statement that this result holds
uniformly relative to $\jump{0}{(3)}$.

\begin{lemma}\label{lem:uniform} There is a computable function $g$ such
that for every
$e$ and $X$, if $\red{e}{X}{ }$ is a $2$--coloring $\color{C}$ of
$\ktup{\omega}{2}$,
then $\red{g(e)}{X^{(3)}}{ }(0)$ is a number $a$ such that $\jump{(H \oplus  
X)}{\prime\prime}
= \red{a}{\jump{X}{\prime\prime}}{ }$ for some infinite set $H$ which is
homogeneous for
$\color{C}$.
\end{lemma}

The lemma is proved by analyzing a relativized version
of the construction in \cite{ChJoSl}, choosing the sets in the
forcing conditions to belong to a fixed Scott set which contains
$X$ and is uniformly low relative to $X$.  The desired uniformity
is established by a quantifier count.

It follows from the lemma that there is a binary function $c \T \jump{0}{(3)}$ 
such that if $A$ is an infinite set with
$A'' = \red{e}{\zeropp}{ }$ and
$\varphi_a$ is a $2$--coloring $\color{C}$ of $\ktup{\omega}{2}$, then there is
a $\color{C}$--homogeneous set $H \subseteq A$ such that $H$ is infinite
and $H'' = \red{c(e,a)}{\zeropp}{ }$.  (This is shown by using the function
$p_A$, which enumerates $A$ in increasing order, to  ``pull back'' the
restriction of $\varphi_a$ to $\ktup{A}{2}$ to an $A$--computable
$2$--coloring of
$\ktup{\omega}{2}$, to which the previous lemma can be applied, and then
taking the image under $p_A$ of the low$_2$ homogeneous set yielded by
the lemma.  In other words, we form a $2$--coloring $\widehat{\color{C}}$ of  
$\ktup{\omega}{2}$ as follows; if $m < n$, define the
$\widehat{\color{C}}$--color of $\{m,n\}$ to be the same as the
$\color{C}$--color of $\{p_A(m), p_A(n)\}$.  Then $\widehat{\color{C}} \T A$  
by a known index, and we apply \lemref{lem:uniform} with $X = A$ to obtain an  
infinite $\widehat{\color{C}}$--homogeneous set $\widehat{H}$ as described
in the lemma.  The set $H = p_A(\widehat{H})$ is then an infinite
$\color{C}$--homogeneous set, and $H \T \widehat{H} \join A$ by a known
index.)   By iterating $c$, one obtains a sequence $A_0  \supseteq A_1
\supseteq \dots$  of infinite sets and a function $h
\T \jump{0}{(3)}$ such that for all $e$, $\jump{A_e}{\prime\prime} =
\red{h(e)}{\zeropp}{ }$
and, if $\varphi_e$ is a $2$--coloring of $\ktup{\omega}{2}$, then $A_e$
is homogeneous for $\varphi_e$.

To complete the proof, it clearly suffices to construct a
$\Delta^0_3$ set $C$ such that $C \subseteq^* A_e$ for all $e$
and $K \not \T C$.
We first show how to construct such a $C$
where we require that $C$ be $\Delta^0_4$ instead of
$\Delta^0_3$.  Then a finite injury modification of the argument
will produce such a $C$ which is $\Delta^0_3$.  We obtain $C$ as
$\cup_e \sigma_e$, where $\sigma_0 \subseteq \sigma_1 \subseteq \dots$ and
each $\sigma_e$ is a binary string.  Let $\sigma_0$ be the empty string.
Given $\sigma_e$, we choose $\sigma_{e+1} \supseteq \sigma_e$ so
as to ensure that $K \neq \{e\}^C$.
Call a string $\tau$ {\em
$e$--acceptable} if $\tau$ extends $\sigma_e$ and $(\forall
k)[\tau(k) = 1 \rightarrow (k \in A_e \vee \sigma_e(k) = 1)]$.

\setcounter{case}{0}
\begin{case}  There exist $e$--split strings $\tau_1$ and $\tau_2$ which
are each $e$--acceptable; i.e., there exist $e$--acceptable strings $\tau_1$  
and $\tau_2$ and an $n$ with $\red{e}{\tau_1}{ }(n)\downarrow \ne
\red{e}{\tau_2}{ }(n)\downarrow$.
\end{case}
Then let $\sigma_{e+1}$ be the string
$\tau$ of least G\"odel number such that $\tau$ is
$e$--acceptable, $\{e\}^\tau$ is incompatible with $K$, and
$\tau(x) = 1$ for some $x \geq |\sigma_e|$.  It is easily seen
that such a string $\tau$ exists, since either $\tau_1$ or
$\tau_2$ has all the requisite properties of $\tau$ except
possibly the last, and then $\tau$ may be chosen as an extension
of $\tau_1$ or $\tau_2$ such that $\tau(x) = 1$ for some $x \in
A_e$ with $x > |\sigma_e|$, and $\tau$ takes the value $0$ on all
other arguments $n > |\tau_i|$ for the appropriate $i$.

\begin{case} Otherwise.
\end{case}
Then let $\sigma_{e+1}$ be the
$e$--acceptable string $\tau$ of least G\"odel number such that
there exists $x \geq |\sigma_e|$ with $\tau(x)=1$.

We now verify that this construction works.  Let $C = \cup_e
\sigma_e$.  Note that if $i \geq e$ then, as $A_i \subseteq A_e$
and $\sigma_i \supseteq \sigma_e$, every $i$--acceptable string
is $e$--acceptable.   It follows that $C \subseteq^* A_e$ for
each $e$ and hence that $C$
is $2$--r--cohesive.  Also, the set of
$e$--acceptable strings is computable from $A_e$, uniformly in
$e$, so the division between Case~1 and Case~2 in the definition
of $\sigma_{e+1}$ is $\Sigma_1^{A_e}$, uniformly in $e$.  The
choice of $\sigma_{e+1}$ in each case can be computed from $A_e
\oplus K$, uniformly in $e$.  As $\jump{A_e}{\prime\prime} \T \zeropp$,
uniformly in
$0^{(3)}$, it follows that $C \T \jump{0}{(3)}$.  Finally, assume
for a contradiction that $K = \red{e}{C}{ }$.  Then Case~2 applies in
the definition of $\sigma_{e+1}$.  It follows that $K \T
A_e$, which is a contradiction since $A_e$ is $\low{2}$.  To see this, note  
that to compute $K(x)$, one need only find any
$e$--acceptable $\tau$ with $\red{e}{\tau}{ }(x) \downarrow$ and then
$K(x) = \red{e}{\tau}{ }(x)$, since there is an $e$--acceptable $\mu$
with $\red{e}{\mu}{ }(x) \downarrow$ and Case~1 does not apply.

The construction of a $2$--r--cohesive set $\Delta^0_3$ set $C$
with $K \not \T C$ is similar, except that the
$\jump{0}{(3)}$--computable functions which occur in the previous
argument are approximated by suitable $\zeropp$--computable
functions.  Each requirement is satisfied essentially as before,
but its satisfaction may need to await the convergence of these
approximations.  We use the same sequence $A_0, A_1, \dots$, but
now let $h$ be a $\jump{0}{(3)}$--computable function with $A_e =
\red{h(e)}{\zeropp}{ }$ for all $e$.  Let $\hat{h}$ be a
$\zeropp$--computable function  with $\lim_s \hat{h}(e,s) = h(e)$ for all
$e$.  To deal with the case distinction in defining
$\sigma_{e+1}$ above, we define for each $e$ a set $B_e$.  Let
\[B_e = \{\sigma \mid (\exists \tau_1, \tau_2 \supseteq \sigma)[
\tau_1, \tau_2 \hbox{ are $e$--split }\]
\[\land \quad
(\forall x)(\forall j \in \{1,2\})[\tau_j(x) = 1 \implies x
\in A_e \vee \sigma(x) = 1]]\}.\]
Note that $B_e$ is c.e.~in $A_e$,
uniformly in $e$.  Since $\jump{A_e}{\prime\prime} \T \zeropp$, uniformly in 
$\jump{0}{(3)}$, it follows that $B_e \T \zeropp$, uniformly in
$\jump{0}{(3)}$.  Let $q$ be a $\jump{0}{(3)}$--computable function with
$\red{q(e)}{\zeropp}{ } = B_e$ for all $e$, and let $\hat{q}$ be a
$\zeropp$--computable function with $\lim_s \hat{q}(e,s) = q(e)$ for
all $e$.

     We now give the construction of $C$, which is carried out
with a $\zeropp$ oracle.  The construction produces a sequence $\sigma_0
\subseteq \sigma_1 \subseteq \dots$  of
strings, and we let $C =
\cup_e \sigma_e$.  It follows that $C \T \zeropp$, so that
$C$ will be $\Delta^0_3$.  The requirement $R_e : K \neq
\red{e}{C}{ }$ is said to be {\em satisfied} at stage $s+1$ if there is
a stage $t \leq s$ such that $R_e$ received attention at stage
$t$, $\hat{h}(i,u) = \hat{h}(i,t)$ for all $i \leq e$ and all $u$
with $t \leq u \leq s+1$, and, finally, $\hat{q}(e,u) =
\hat{q}(e,t)$ for all $u$ with $t \leq u \leq s+1$.  The requirement $R_e$
is said to {\em require attention\/} at stage $s+1$ if it is not satisfied at  
stage $s+1$.

\noindent {\em Stage $0$.}  Let $\sigma_0 = \emptyset$.

\noindent {\em Stage $s+1$.}  Let $\sigma_s$ be given.  Let $e$ be the least 
number such that $R_e$ requires attention.  (Such a number exists
because there are only finitely $e$ such that $R_e$ does not
require attention.)  Let $t$ be the least number such that $t
\geq s$ and at least one of the following conditions holds:

{\renewcommand{\labelenumi}{(\roman{enumi})}
\begin{enumerate}
\item $\hat{h}(i,t) \neq \hat{h}(i,t+1)$ for some $i \leq e$
\item $\hat{q}(e,t) \neq \hat{q}(e,t+1)$
\item there exist $e$--split strings $\tau_1$ and $\tau_2$, each
extending $\sigma_s$, such that
\[(\forall x)(\forall j \in \{1,2\})[\tau_j (x) = 1 \implies x < |\sigma_s|  
\vee (\forall i \leq e)[\red{\hat{h}(i,t)}{\zeropp}{t}(x) = 1]]\]
and
\[(\forall j \in \{1,2\})(\exists x \geq |\sigma_s|)[\tau_j(x) = 1]\]
\item $\red{\hat{q}(e,t)}{0''}{t} (\sigma_s) = 0$.
\end{enumerate}
}
We first note that such a number $t$ must exist.  Suppose for a
contradiction that no $t \geq s$ satisfying at least one of (i) -- (iv)
exists.  Then,
since no $t \geq s$ satisfies (ii) or (iv), it follows that
$\sigma_s \in B_e$.  Consider strings $\mu_1$ and $\mu_2$ which
witness that $\sigma_s \in B_e$.  Since no $t \geq s$
satisfies (i), it is easily seen that $\mu_1$ and  $\mu_2$ satisfy the
first conjunct of (iii).  Now, given $j \in \{1,2\}$, take $x \ge |\mu_j|$
such that $x \in \cap_{i \leq e} A_i$, and define
$\tau_j$ to be an extension of $\mu_j$ such that $x$ is the
unique number $y \geq |\mu_j|$ with $\tau_j(y) = 1$.  It is easily
seen that (iii) holds for $\tau_1$ and $\tau_2$ for all sufficiently large $t$.  

Now we give the definition of $\sigma_{s+1}$.  If (iii) holds,
choose $\sigma_{s+1} \in \{\tau_1, \tau_2\}$ so that
$\{e\}^{\sigma_{s+1}}$ is incompatible with $K$.  Otherwise, search
for $x$ and $t$ such that $t \geq s$, $x \geq |\sigma_s|$, and $(\forall
i \leq e)[\red{h(i,t)}{\zeropp}{t}(x) = 1]$.  Such an $x$ and $t$
exist because $\cap_{i \le e} A_i$ is infinite.   Fix the first
such pair $(x,t)$, and let $\sigma_{s+1}$ be an extension of
$\sigma_s$ such that $x$ is the unique $y \geq |\sigma_s|$ with
$\sigma_{s+1}(y) = 1 $.  If (iii) or (iv) holds,  say that $R_e$
{\em received attention} at stage $s+1$.  This completes the
construction.

To verify that the construction succeeds, first note that it is
$0''$--computable, and so $C$ is $\Delta^0_3$.  In addition, $C$
is infinite because for each $s$ there exists $x \geq |\sigma_s|$
with $\sigma_{s+1}(x) = 1$.   Each $R_e$  requires attention
at most finitely often, since $\lim_s \hat{h}(i,s)$ exists for
each $i \leq e$ and $\lim_s \hat{q}(e,s)$ exists, and $R_e$
receives attention at most once after the least stage $s_0$ such that
$\hat{h}(i,t) = \hat{h}(i,s_0)$ and $\hat{q}(e,t) = \hat{q}(e,s_0)$ for all  
$i \le e$ and all $t \ge s_0$. It
now follows from the construction that $C \subseteq^* A_e$ for
each $e$.  It remains to show that $K \neq \red{e}{C}{ }$ for each $e$.
By construction there is a unique stage $s_1
\geq  s_0$ such that $e$ is the least number such that $R_e$
requires attention at $s_1$. Since $s_1 \geq s_0$, neither (i)
nor (ii) can apply.  If (iii) applies, it is clear that $R_e$ is
satisfied.  Thus we may assume that (iv) holds, from which it
follows that $\sigma_s \notin B_e$.  From this we can argue (as
in the construction for $C \in \Delta^0_4$) that $K \T A_e$,
which is a contradiction because $A_e$ is low$_2$.

\end{proof}

It is not known whether for each noncomputable set $D$ there is a
$\Delta^0_3$ $2$--r--cohesive $A$ with $D \not \T A$, which would
be the natural common generalization of  \thmref{thm:2rcozp} and
\thmref{thm:d32rco}.  However the proof of \thmref{thm:d32rco}
shows this holds for $D$ non--low$_2$, since $D$ may be assumed
without loss of generality to be $\Delta^0_3$ by the existence of
a $\Delta^0_3$ $2$--r--cohesive set. It is also open whether the
arithmetical complexity in \thmref{thm:d32rco} can be  improved
to $\Pi^0_2$, which would be best possible.  Relativizing this
result to $\zerop$ yields the following.

\begin{cor}  There exists a $\Delta^0_4$ $2$--cohesive set $A$ with $\zeropp  
\not\T A$.
\end{cor}

\section{$n$-cohesive and $n$-r-cohesive sets}\label{sec:general}

Here we study arithmetical definability and degrees of
$n$--cohesive and $n$--r--cohesive sets for arbitrary $n$.  Some results
known for the case $n \leq 2$ generalize to arbitrary $n$, whereas
others fail or are open for arbitrary $n$.  We showed in
\thmref{thm:pi22cozpp} that the analogue of the existence of an
incomplete co--maximal set fails for $n=2$.  On the other hand the
following theorem shows that the analogue of the existence of a
cohesive set $A \T K$ holds for all $n$.

\begin{thm}\label{thm:arithnco}  For all $n \geq 1$, there exists
an $n$--cohesive set $B \T \jump{0}{(n)}$.
\end{thm}

\begin{proof}  The following relativized version of the result is
proved by induction on $n$.  For all $n \geq 1$, all $X \subseteq
\omega$, and all infinite $Y \leq_T X$, there is an $n$--cohesive set
$B \subseteq Y$, relative to $X$, such that $B \oplus X \T
\jump{X}{(n)}$.  The result holds for the base case $n = 1$ by
relativizing the result that every infinite computable set
contains a cohesive set $B \T 0'$.  Now assume the
result for $n$ in order to prove it for $n+1$.  To simplify the
notation, we actually prove it for $n+1$ in the unrelativized
case where $X = \emptyset$ and $Y = \omega$.

Fix a $\low{2}$ cohesive set $A$, and recall that the restriction
of any c.e.~coloring of $\ktup{\omega}{n+1}$ to $\ktup{A}{n+1}$
is stable in the sense of \defref{def:stable},  by the obvious
generalization of \lemref{lem:cohstable} to $(n+1)$--tuples.   By
the inductive hypothesis (with $X = A'$ and $Y = A$), take $B
\subseteq A$ such that $B$ is $n$--cohesive relative to
$\jump{A}{\prime}$ and $B \oplus A' \T \jump{(A^\prime)}{(n)} =
\jump{A}{(n+1)} \T \jump{0}{(n+1)}$. Let $W_e$ be a c.e.~subset
of $\ktup{\omega}{n+1}$, viewed as a c.e.~$2$--coloring of
$\ktup{\omega}{n+1}$.  When restricted to $\ktup{A}{n+1}$, this
$2$--coloring is stable and so induces a $\Sigma^0_2(A)$
$2$--coloring $\color{P}_e$ of $\ktup{A}{n}$.  Specifically, let
$\color{P}_e$ be the set of all $D \in \ktup{A}{n}$ such that there
are only finitely many $a \in A$ with $D \cup \{a\} \in W_e$.  Then,
by stability, if $D \in \ktup{A}{n}$ and $D \notin
\color{P}_e$, there are only finitely many $a \in A$ such that $D
\cup \{a\} \notin W_e$.  By the choice of $B$, $B$ is almost
homogeneous for each partition $\color{P}_e$.

To complete the proof, it suffices to show that there is a set $C
\subseteq B$ such that $C$ is $(n+1)$--cohesive and $C \T
B \oplus \jump{A}{\prime\prime} \T \jump{0}{(n+1)}$.  Let $C_s$ be the first  
$s$ elements of
$C$ in natural order.  Suppose that we have already defined $C_s$.
Let $c_s$ be the least number $z > \max C_s$ such that $z \in B$ and
for every $D \in \ktup{C_s}{n}$ and every $e \leq s$, $D \in
\color{P}_e$ if and only if $D \cup \{z\} \notin W_e$.  Such a
number $z$ exists because all sufficiently large elements of $A$
(hence of $B$) have the desired properties for $z$, and
furthermore $c_s$ can be found by a $B \oplus A''$--effective
search.  Let $C_{s+1} = C_s \cup \{c_s\}$. It is clear that $C
\subseteq B$ and that $C \T  B \oplus \jump{A}{\prime\prime}$.

It remains to check that $C$ is $(n+1)$--cohesive.  Let a
c.e.~$2$--coloring $W_e$ of $\ktup{\omega}{n+1}$ be given, and recall that
$\color{P}_e$ is the induced $2$--coloring of $\ktup{A}{n}$.  Since $B$ is
almost homogeneous for $\color{P}_e$ and $B \subseteq A$, we may
choose a finite set $F$ such that $B - F$ is homogeneous for
$\color{P}_e$. Let $G = C - (F \cup C_e)$.  We claim that $G$ is
homogeneous for $W_e$, which completes the proof.  Suppose first that
$\ktup{B - F}{n} \subseteq \color{P}_e$.  Fix $D \in \ktup{G}{n}$, so $D \in
\color{P}_e$.  Let $c_s \in C$ with $c_s > \max D$, so that $D
\subseteq C_s$.  Then $s > e$, since otherwise $D \subseteq C_s
\subseteq C_e$, in contradiction to $D$ being a nonempty set disjoint
from $C_e$.  It then follows from the choice of $c_s$ that $D \cup
\{c_s\} \notin W_e$.  As any $(n+1)$--element subset of $G$ can
be written in the form $D \cup \{c_s\}$ with $D \in \ktup{G}{n}$ and
$c_s > \max D$, this argument shows that $\ktup{G}{n+1}$ is disjoint
from $W_e$ on the assumption that $\ktup{B - F}{n} \subseteq
\color{P}_e$.  An entirely analogous argument shows that
$\ktup{G}{n+1}$ is contained in $W_e$ if $\ktup{B-F}{n}$ is disjoint from
$\color{P}_e$.
\end{proof}

We don't know for which $n \ge 2$ there is an
$n$--cohesive set $B$ with $B <_T \jump{0}{(n)}$.  If it could be
shown in a relativizable fashion that there is such a set $B$ for
$n = 2$, then there is such a set $B$ for each $n \geq 2$.  This
is proved by essentially the same inductive argument as is used
to prove \thmref{thm:arithnco}, together with the observation
that the set $B$ chosen in the inductive step satisfies $0'' \T
B$ since it is $2$--cohesive relative to $\zerop$.  From this it
follows, in the notation of the proof of \thmref{thm:arithnco}, that
$C \T B \oplus A'' \T B$, where $B <_T \jump{A}{(n+1)}$ by
inductive hypothesis.

As an immediate corollary to \thmref{thm:arithnco}, we have the following.

\begin{cor}  Let $n \ge 1$.  Every c.e.~$2$--coloring of $\ktup{\omega}{n}$   
has an infinite $\Delta^0_{n+1}$
homogeneous set.
\end{cor}
As Jockusch's result \thmref{thm:basicjockusch}, when
relativized to $\zerop$, only  guarantees an infinite $\Pi^0_{n+1}$
homogeneous set, this corollary provides  new information about
infinite homogeneous sets for c.e.~$2$--colorings of
$\ktup{\omega}{n}$, $n \ge 3$.  Recall that, although every
c.e.~$2$--coloring of $\ktup{\omega}{1}$ has an infinite
$\Pi^0_1$  homogeneous set and every c.e.~$2$--coloring of
$\ktup{\omega}{2}$ has an  infinite $\Pi^0_2$ homogeneous set, it
is unknown whether every  c.e.~$2$--coloring of
$\ktup{\omega}{n}$, $n \ge 3$, has an infinite  $\Pi^0_n$
homogeneous set.

We next give some results relating the notions of $n$--cohesive and
$n$--r--cohesive sets to relativizations of the
usual notions of cohesive and r--cohesive sets.  In  what follows,
a $\jump{0}{(n-1)}$--cohesive set (respectively
$\jump{0}{(n-1)}$--r--cohesive set) is a set which is
cohesive  (respectively r--cohesive) relative to
$\jump{0}{(n-1)}$.

\begin{thm}\label{thm:kcoco}  Let $n \ge 1$.  If $A$ is $n$--cohesive, then  
$A$ is $\jump{0}{(n-1)}$--cohesive.
\end{thm}

\begin{proof}  We first prove the following lemma.

\begin{lemma}\label{lem:relcoh}  Let $A$ be any set and let $n \ge 1$.  If a  
set $B$ is
$(n+1)$--cohesive relative to $A$, then $B$ is $n$--cohesive relative to
$\jump{A}{\prime}$.
\end{lemma}

\begin{proof}  Assume that $B$ is $(n+1)$--cohesive relative to $A$.  Let
$\color{C}$ be a $2$--coloring of $\ktup{\omega}{n}$, given as a subset of
$\ktup{\omega}{n}$, which is c.e.~in $\jump{A}{\prime}$; i.e.,
$\Sigma^0_2(A)$.  Then $\color{C} \m \fin{A} = \{e \mid W_e^A \mbox{ is
finite}\}$, say via a computable function $f$; i.e.,
\[\{x_1, \dots, x_n \}\in \color{C} \iff f(x_1, \dots, x_n) \in \fin{A}.\]

Let $\color{C}_s \subseteq \ktup{\omega}{n}$ be the $A$--computable
$2$--coloring of $\ktup{\omega}{n}$ defined by
\[\color{C}_s = \{\{x_1, \dots, x_n \} \mid x_1 < \dots < x_n < s \land
W_{f(x_1, \dots, x_n),s}^A = W_{f(x_1, \dots, x_n),s+1}^A\}.\]
Then $\color{C}$ can be computably approximated, relative to $A$, by
$\{\color{C}_s\}_{s
\in \omega}$; i.e., for all $x_1 < x_2 < \dots < x_n$,
\[\{x_1, \dots, x_n\} \in \color{C} \iff (\exists s_0)(\forall s \ge
s_0)[\{x_1, \dots, x_n\} \in \color{C}_s].\]

Next, define a $2$--coloring $\color{Q} \subseteq \ktup{\omega}{n+1}$ as
follows:  if $x_1 < x_2 < \dots < x_n < x_{n+1}$, then
\[\{x_1, \dots, x_n, x_{n+1}\} \in \color{Q} \iff (\exists s >
x_{n+1})[\{x_1, \dots, x_n\} \notin \color{C}_s].\]
Then $\color{Q}$ is c.e.~in $A$, and since $B$ is $(n+1)$--cohesive relative  
to $A$, there exists a finite set $F$ such that $B - F$ is an infinite
$\color{Q}$--homogeneous set.  We show that $B-F$ is a $\color{C}$--homogeneous  
set.

\setcounter{case}{0}
\begin{case} $\ktup{B-F}{n+1} \subseteq \color{Q}$.
\end{case}

We show that $\ktup{B-F}{n} \subseteq \overline{\color{C}}$.  Let $x_1 <
\dots < x_n$ be elements of $B-F$.  Given $y > x_n$ with $y \in B-F$, $\{x_1,  
\dots, x_n, y\} \in \color{Q}$, which implies that there is $s > y$ such
that $\{x_1, \dots, x_n\} \notin  \color{C}_s$.  Since $B-F$ is infinite, we  
have $\{x_1, \dots x_n\} \notin \color{C}_s$ for infinitely many $s$, and
hence $\{x_1, \dots, x_n\} \notin \color{C}$.  Hence $\ktup{B-F}{n} \subseteq  
\overline{\color{C}}$.

\begin{case}  $\ktup{B-F}{n+1} \subseteq \overline{\color{Q}}$.
\end{case}

We show that $\ktup{B-F}{n} \subseteq \color{C}$.  Let $x_1 < \dots < x_n$
be elements of $B-F$.  Given $y > x_n$ with $y \in B-F$, $\{x_1, \dots,
x_n,y\} \notin \color{Q}$, which implies that for all $s > y$, $\{x_1, \dots,  
x_n\}\in \color{C}_s$.  Hence $\{x_1, \dots, x_n\}\in \color{C}$, and so
$\ktup{B-F}{n} \subseteq \color{C}$.

Hence $B$ is $n$--cohesive relative to $\jump{A}{\prime}$.
\end{proof}

We now prove the theorem by induction on $n$.  If $n = 1$, then we have
already noted that any $1$--cohesive set is cohesive.  So, assume the result  
for $n$ and prove it for $n+1$.  Assume that $B$ is $(n+1)$--cohesive.  By the  
lemma, $B$ is $n$--cohesive relative to $\zerop$.  By the inductive
hypothesis, relativized to $\zerop$, $B$ is cohesive relative to
$\jump{0}{(n)}$.
\end{proof}

Note that the converse of \thmref{thm:kcoco} is not true.  Relativizing the  
result that there exists an incomplete maximal set (i.e., there exists a
$\Pi_1^0$ cohesive set $A$ such that $A <_T \zerop$) to $\zerop$, we see that  
there exists a $\Pi_2^0$ $\zerop$--cohesive set $A$ such that $A <_T
\zeropp$.  Since $\Pi_2^0$ $2$--cohesive sets are of degree $\zeropp$ by
\thmref{thm:pi22cozpp}, it follows immediately that there exists a
($\Pi_2^0$) $\zerop$--cohesive set which is not $2$--cohesive.

A result similar to \thmref{thm:kcoco} can be proved for $n$--r--cohesive sets.

\begin{thm}\label{thm:krcorco}  Let $n \ge 1$. If $A$ is $n$--r--cohesive,
then $A$ is $\jump{0}{(n-1)}$--r--cohesive.
\end{thm}

The proof of this theorem is analogous to that of
\thmref{thm:kcoco} but uses the Limit Lemma.

We call a degree $\bf a$ $n$--{\em cohesive} if there is an
$n$--cohesive set of degree $\bf a$, and the $n$--r--cohesive degrees
are defined analogously.  Our next theorem will give some information
on $n$--cohesive and $n$--r--cohesive degrees and their jumps for $n
\geq 2$, but we first recall some results for the case $n = 1$.  To
state these we need a definition.

\begin{definition} Let $\bf a$ and $\bf b$ be
degrees.  Then ${\bf a << b}$ means that any partial $\bf
a$--computable $\{0,1\}$--valued function can be extended to a total
$\bf b$--computable function.  (This notation is due to S. Simpson
\cite{Si}, page 648, who used a differently stated but equivalent definition.)
\end{definition}

It was shown by Jockusch and Stephan \cite{JoSt1}, Corollary 2.4, that
the cohesive degrees coincide with the r--cohesive degrees.  Also, it
was shown in \cite{JoSt1}, Theorem 2.2(ii) (see also \cite{JoSt2}),
that the jumps of the r--cohesive degrees are precisely the degrees
${\bf c >> 0'}$.  The following easy lemma is a consequence of this
latter result.

\begin{lemma}\label{lem:jcoh}  Assume that $D \not \leq_T 0'$.  Then
there is an r--cohesive set $A$ with $D \not \leq_T A'$.
\end{lemma}

\begin{proof}
First note that there is a $\Pi^0_1$ class $P \subseteq 2^\omega$ such
that the degrees of elements of $P$ are exactly the degrees ${\bf a >>
0}$;  namely, we let
\[P = \{f \in 2^\omega \mid (\forall e)(\forall i)[\varphi_e(i) \downarrow
\leq 1 \implies f(\langle e,i \rangle) = \varphi_e(i)] \}.\]
Hence, by a result \cite{JoSo} (Theorem 2.5) of Jockusch and Soare on
$\Pi^0_1$ classes, for any degree ${\bf d > 0}$ there is a degree
${\bf a >> 0}$ with ${\bf a \not \geq d}$.  Relativizing this result
to $\bf 0'$ and applying the result \cite{JoSt1} (Theorem 2.2(ii))
that the jumps of the r--cohesive degrees are precisely the degrees
${\bf c >> 0'}$, yields the lemma.
\end{proof}

Our next theorem shows that the $n$--cohesive and $n$--r--cohesive
degrees behave quite differently for $n \geq 2$ than for $n = 1$ as
described above.  The special case where $n = 2$ has been treated in
earlier sections of this paper.

\begin{thm} \label{thm:hier}  Assume that $n \geq 2$.
\begin{enumerate}
\item If $A$ is $n$--r--cohesive, then $A \not \T \jump{0}{(n-1)}$.
\item If $A$ is $n$--r--cohesive, then $\jump{0}{(n-2)}
\T A$. Conversely, if $D \not\T \jump{0}{(n-2)}$ then there exists
an $n$--r--cohesive set $A$ such that $D \not\T A$.
\item  If $A$ is $n$--cohesive, then $\jump{0}{(n-1)}
\T A$.  Conversely, if $D \not\T \jump{0}{(n-1)}$ then there exists
an $n$--cohesive set $A$ such that  $D \not\T A$.
\item  ${\bf \jump{0}{(n+1)}}$ is the least degree among all jumps of
$n$--cohesive degrees.
\item  If $A$ is $n$--r--cohesive, then $\jump{0}{(n)} <_T
\jump{A}{\prime}$.

\end{enumerate}
\end{thm}

\begin{proof}
The first statement follows immediately from the second part of
\thmref{thm:basicjockusch}.  The first part of the second statement
follows immediately from a result of Jockusch \cite{Jo1}, Lemma 5.9,
which asserts that for each $n \geq 2$ there is a computable
$2$--coloring $\color{C}$ of $\ktup{\omega}{n}$ such that
$\jump{0}{(n-2)}\leq_T A$ for every infinite $\color{C}$--homogeneous
set $A$.

The second part of the second statement is proved by induction on $n
\geq 2$ in the following relativized form: for any set $X \subseteq
\omega$, any infinite set $Y \leq_T X$, and any set $D \not \leq_T
\jump{X}{(n-2)}$, there is a set $A \subseteq Y$ such that $A$ is
$n$--r--cohesive relative to $X$ and $D \not\T A \join X$.  For $n =
2$, this follows from the relativization of \thmref{thm:2rcozp}.

The inductive step is similar to the proof of \thmref{thm:arithnco}.
Assume our result is true for $n$, where $n \geq 2$.  We prove it for
$n+1$ assuming, for notational simplicity, that $X = \emptyset$ and $Y
= \omega$.  Let $D \not \leq_T \jump{0}{(n-1)}$ be given.  We must
construct an $(n+1)$--r--cohesive set $A$ with $D \not
\leq_T A$.  One obstacle to constructing $(n+1)$--r--cohesive sets is
that the computable sets are not uniformly computable.  To overcome
this obstacle, let $L$ be a set of low degree such that there is a
sequence of sets $R_0, R_1, \dots$ with $\bigoplus \{R_i \mid i \in
\omega\} \T L$ and such that every computable subset of
$\ktup{\omega}{n+1}$ (i.e., every computable $2$--coloring of
$\ktup{\omega}{n+1}$) occurs as $R_i$ for some $i$.  (The existence of
such an $L$ follows from an application of the low basis theorem
\cite{JoSo} (Theorem 2.1) to the $\Pi^0_1$ class $P$ constructed in
the proof of
\lemref{lem:jcoh}.   Namely, let $f$ be an element of $P$ of low
degree and let $R_e = \{i \mid f(\langle e,i\rangle) = 1\}$.)
Let $U$ be a set which is r--cohesive relative to $L$  with
$D \not \leq_T  \jump{(U \oplus L)}{(n-1)}$.  If $n = 2$, the
existence of such a set $U$ follows from \lemref{lem:jcoh}
relativized to $L$ and the fact that $D \not \T L'
\equiv_T 0'$.   If $n \geq 3$, choose $U$ to be r--cohesive
relative to $L$ and low$_2$ relative to $L$.  It follows that $U
\oplus L$ is low$_2$, since $L$ is low$_2$.  Then $D \not
\leq_T \jump{0}{(n-1)} \equiv_T \jump{(U \oplus L)}{(n-1)}$,
since $n \geq 3$.

By the inductive hypothesis (with $X = (U \oplus L)'$ and $Y = U$),
there is a set $B \subseteq U$ such that $B$ is $n$--r--cohesive
relative to $(U \oplus L)'$ and $B
\oplus (U \oplus L)' \not \geq_T D$.  Note that the restriction of
each $R_e$ to $\ktup{U}{n+1}$  is stable, since $U$ is
r--cohesive,
relative to $L$.  Hence each $R_e$ induces a $2$--coloring $S_e$ of
$\ktup{U}{n}$
(thought of as a subset of  $\ktup{U}{n}$)
such
that for all  $D \in \ktup{U}{n}$
and all sufficiently large $t \in U$,
$D \in S_e$ if and only if $D \cup \{t\} \in R_e$.  Note that, by the
Limit Lemma, $S_e \leq_T (U \oplus L)'$, uniformly in $e$.

We now use the technique of \thmref{thm:arithnco} to construct an
$(n+1)$--r--cohesive set $C \subseteq B$ such that $C
\leq_T B \oplus (U \oplus L)'$.  Let $C_s$ be the first $s$ elements of $C$ in
natural order.  Suppose that we have already defined $C_s$.  Let $c_s$
be the least number $z > \max C_s$ such that $z \in B$ and for every
$D \in \ktup{C_s}{n}$ and every $e \leq s$, $D \in S_e$ if and only if
$D \cup \{z\} \in R_e$.  Such a number $z$ exists because all
sufficiently large elements of $U$ (hence of $B$) have the desired
properties for $z$, and furthermore $c_s$ can be found by a $B
\oplus (U \oplus L)'$--effective search.  Let $C_{s+1} = C_s \cup
\{c_s\}$. It is clear that $C \subseteq B$ and that $C \T  B \oplus
(U \oplus L)'$.  This completes the construction of $C$.

We claim now that $C$ is the desired $(n+1)$--r--cohesive set with $D
\not \leq_T C$.  The proof that $C$ is $(n+1)$--r--cohesive is entirely
analogous to the corresponding verification in the proof of
\thmref{thm:arithnco}, so we omit it.  Assume now for a contradiction
that $D \leq_T C$.  Then $D \leq_T C \leq_T B \oplus (U \oplus L)'$,
which contradicts our choice of $B$.  This completes our proof of the
second statement in the second part of the theorem.

The first part of the third statement is proved in the following
relativized form by induction on $n$: if $A$ is $n$--cohesive relative
to $X$, then $\jump{X}{(n-1)} \T A \join X$.  This is obvious for
$n=1$; for $n=2$ it is proved by relativizing
\propref{prop:2coabovezp}.  The inductive step uses
\lemref{lem:relcoh} and the case $n = 2$.  To prove the second
statement in the third part of the theorem, assume that $D \not
\leq_T \jump{0}{(n-1)}$.  By the second statement in the second part,
relative to $0'$, there is a set $A$ which is $n$--r--cohesive
relative to $0'$ with $D \not \leq_T A$.  Clearly $A$ is
$n$--cohesive, so the proof of the third part of the theorem is
complete.  (The second statement of this part could alternatively be
proved directly by induction on $n$ along the lines of
\thmref{thm:arithnco}.)

To prove the fourth part of the theorem, we first show that for any
$n$--cohesive set $A$, $\jump{0}{(n+1)} \T
\jump{A}{\prime}$. This is proved in relativized form
by induction on $n \ge 2$.  The base case $n=2$ is obtained by
relativizing \propref{prop:high2co}.  The inductive step uses
\lemref{lem:relcoh} and the case $n=2$.  Thus, ${\bf a' \geq
\jump{0}{(n+1)}}$ for any $n$--cohesive degree $\bf a$.  It remains to
show there is an $n$--cohesive degree $\bf a$ with ${\bf a' =
\jump{0}{(n+1)}}$.  By \thmref{thm:arithnco} there is an $n$--cohesive
degree ${\bf a \leq \jump{0}{(n)}}$.  For this $\bf a$, we have ${\bf a' \leq
\jump{0}{(n+1)}}$ since ${\bf a \leq \jump{0}{(n)}}$, and ${\bf a'
\geq \jump{0}{(n+1)}}$ by what is proved just above.  Thus ${\bf a' =
\jump{0}{(n+1)}}$ as needed.  (In fact, as we will remark later, we
could simply choose ${\bf a = \jump{0}{(n)}}$.)

We prove the fifth and final part of the theorem by induction on
$n \ge 2$ in the following relativized form:  if $A$ is
$n$--r--cohesive  relative to $X$, then $\jump{X}{(n)} <_T
\jump{(A \join X)}{\prime}$.  The  base case $n=2$ is obtained
from \thmref{thm:high2rco}, relativized to $X$.   The inductive
step, with the usual simplifying assumption that $X =
\emptyset$, uses the analogue of \lemref{lem:relcoh} and
the fact that  if $A$ is $(n+1)$--r--cohesive, $n \ge 2$, then
$\zerop \T A$.

\end{proof}

We now consider implications among various forms of cohesiveness, both
for sets and for degrees.  It is obvious that every $n$--cohesive set
is $n$--r--cohesive, and it is easy to check that if $n \geq k$, then
every $n$--cohesive set is $k$--cohesive, and every $n$--r--cohesive
set is $k$--r--cohesive.  Furthermore, each $(n+1)$--r--cohesive set
is $n$--r--cohesive relative to $0'$ by the analogue to
\lemref{lem:relcoh} and hence
is $n$--cohesive.  These simple remarks imply that our hierarchy of
notions of cohesiveness is linearly ordered:

$1$--r--cohesive, $1$--cohesive, $2$--r--cohesive, $2$--cohesive,
$3$--r--cohesive, \dots

Of course, the notions become stronger as one goes from left to
right.  The first inclusion ($1$--cohesive implies $1$--r--cohesive) is
proper for sets (by \cite{So1}, page 191) but not for degrees (by
\cite{JoSt1}, Corollary 2.4).  The following corollary shows that all
the remaining inclusions are proper, even for degrees.

\begin{cor}\label{cor:stronghier}  Assume $n \geq 2$.
\begin{enumerate}
\item There is an $n$--r--cohesive degree which is not $n$--cohesive.
\item  There is an $(n-1)$--cohesive degree which is not
$n$--r--cohesive.
\end{enumerate}
\end{cor}

\begin{proof}
The first part of the corollary follows at once from the second and
third parts of \thmref{thm:hier}  because the second part of
\thmref{thm:hier} implies that there exists an $n$--r--cohesive set $A$ such  
that $\jump{0}{(n-1)} \not\T A$. To prove the second part  of the
corollary, consider an $(n-1)$--cohesive set $A \leq_T
\jump{0}{(n-1)}$, which exists by \thmref{thm:arithnco}.  Then
$\deg(A)$ is not  $n$--r--cohesive  by the first part of
\thmref{thm:hier}. \end{proof}

The following result extends \thmref{thm:d32rco} and in particular
implies that for each $n \geq 2$ there is an $n$--r--cohesive degree
${\bf a < \jump{0}{(n)}}$.

\begin{thm}
For each $n \geq 2$ there exists an $n$--r--cohesive set $A$ such that
$A \T \jump{0}{(n)}$ and $\jump{0}{(n-1)} \not \T A$.
\end{thm}

\begin{proof}  This result is proved in relativized form by induction
on $n \geq 2$.  The base step where $n = 2$ is the relativized
version of \thmref{thm:d32rco}.  The induction step is almost
identical to the induction step in the proof the second part of the
second statement in \thmref{thm:hier} (but does not need a separate
treatment of the case when $n=2$), so we omit the details.
\end{proof}

We note that for each $n \geq 1$ the $n$--cohesive and
$n$--r--cohesive degrees are closed upwards by a result of Jockusch
\cite{Jo2} (Theorem~1) and \thmref{thm:arithnco} of the current paper.
Thus, by \thmref{thm:arithnco}, for each $n \geq 1$ the degree ${\bf
\jump{0}{(n)}}$ is $n$--cohesive.

Furthermore, our results characterize the jumps of the
$n$--cohesive degrees, for $n \ge 2$; namely, for $n\ge 2$, a
degree $\td{a}$ is the jump of an $n$--cohesive degree if and
only if $\td{a} \ge {\bf\jump{0}{(n+1)}}$. This follows
immediately from \thmref{thm:hier}, part (4), and the upward
closure of the jumps of the $n$--cohesive degrees.   The latter upward
closure result
follows immediately from the upward closure of the $n$--cohesive
degrees and the relativized version of the Friedberg
Completeness Criterion (see \cite{So1}, Theorem~VI.3.2).

It remains open to find a characterization of the $n$--cohesive
or $n$--r--cohesive degrees for any $n \geq 1$,  or to
find a characterization of the jumps of the $n$--r--cohesive
degrees, for $n \ge 2$.  It is consistent with \thmref{thm:hier}
(but seems unlikely) that for $n \geq 2$ the $n$--cohesive
degrees are precisely the degrees $\bf a$ such that ${\bf a \geq
\jump{0}{(n-1)}}$ and ${\bf a' \geq \jump{0}{(n+1)}}$.

% References

\end{document}